\documentclass[11pt, reqno]{amsart}
\usepackage{amsfonts, amsmath, amscd,amssymb}
 \usepackage{fullpage}
 
\usepackage{amssymb, amsbsy, amsthm,  amstext, amsopn, verbatim,multicol}
\usepackage[mathscr]{eucal}
\usepackage[all]{xy}
\usepackage{mathrsfs}
\usepackage{latexsym}

\makeatletter
\renewcommand{\subsection}{\@startsection{subsection}{2}{0pt}{-3ex plus -1ex minus -0.2ex}{-2mm 
plus -0pt minus -2pt}{\normalfont\bfseries}} \makeatother

\numberwithin{equation}{subsection}
\newtheorem*{theorem}{Theorem}
\newtheorem*{proposition}{Proposition}
\newtheorem*{lemma}{Lemma}

\newtheorem*{corollary}{Corollary}

\theoremstyle{remark}
    
\makeatletter
\newcommand{\cmpl}[1]{%
    \sbox\z@{$#1$}%
    \dimen@=\wd\z@
    \advance \dimen@ -\strip@pt\fontdimen\@ne\textfont\@ne \ht\z@
    \setbox\tw@=\hb@xt@\dimen@{}%
    \ht\tw@=\ht\z@ \dp\tw@=\dp\z@
    \box\z@
    \llap{$\overline{\box\tw@}$}%
}
\makeatother

\newcommand{\Lmod}[1]{#1\text{-}{\mathsf{mod}}}

\newcommand{\Lmof}[1]{#1\text{-}{\mathsf{grmod}}}

\DeclareMathOperator{\Ext}{\mathrm{Ext}}
\DeclareMathOperator{\Hom}{\mathrm{Hom}}
\DeclareMathOperator{\End}{\mathrm{End}}
\DeclareMathOperator{\im}{\mathrm{Im}}
\DeclareMathOperator{\qgr}{{\mathbb F}}

\DeclareMathOperator{\Rad}{\mathfrak{R}}   
\DeclareMathOperator{\Rd}{\cmpl{\mathfrak{R}}}
\DeclareMathOperator{\Rdd}{\mathrm{Rad}}

\DeclareMathOperator{\gr}{\mathrm{gr}}

\DeclareMathOperator{\supp}{{\mathsf{Supp}}}

\DeclareMathOperator{\Tr}{\mathrm{Tr}}

\DeclareMathOperator{\gk}{\mathrm{GKdim}}

\newcommand{\erem}{\hfill$\lozenge$\end{remark}}
\newcommand{\beq}{\begin{equation}\label}
\newcommand{\eeq}{\end{equation}}
\newcommand{\en}{\enspace}

\renewcommand{\o}{\otimes}
\newcommand{\f}[1]{\mathfrak{#1}}
\newcommand{\scr}[1]{\mathscr{#1}}

\newcommand{\iso}{{\ \stackrel{_\sim}{\longrightarrow}\ }}
\def\eu{{\boldsymbol{1}}}

\def\ccirc{{{}_{^{\,^\circ}}}}
\DeclareMathOperator{\Proj}{\mathrm{Proj}}

\newcommand{\sseq}{\subseteq}
\newcommand{\sminus}{\smallsetminus}
\renewcommand{\mid}{\enspace\big|\enspace}
\newcommand{\too}{\,\,\longrightarrow\,\,}

\newcommand{\mto}{\longmapsto}

\newcommand{\ontoo}{\, -\!\!\!\twoheadrightarrow\,}
\newcommand{\cd}{\!\cdot\!}

\newcommand{\g}{{\mathfrak{g}}}

\newcommand{\BH}{{\mathbb H}}
\newcommand{\cyc}{^{\operatorname{cyc}}}

\newcommand{\M}{{\mathcal{M}}}
\newcommand{\FM}{{\mathfrak{M}}}

\newcommand{\SC}{{\scr S}}

\newcommand{\oo}{{\mathcal O}}

\renewcommand{\part}{{\f P}}

\def\C{{\mathbb{C}}}
\def\V{{\f G}}

\def\BS{{\mathbb S}}

\def\gl{{\mathfrak{g}\mathfrak{l}}}
\def\sl{{{\mathfrak{s}\mathfrak{l}}}}

\newcommand{\svol}{\mathbf{s}} 
\newcommand{\PPP}[2]{{}_{#1}{\mathsf{P}}_{#2}}
\newcommand{\QQQ}[2]{{}_{#1}{\mathsf{Q}}_{#2}}

\newcommand{\DDD}[2]{{{}_{#1}{{\mathfrak D}}_{#2}}}

\newcommand{\hccp}{{\mathsf{H}}'_c}
\newcommand{\ehe}{{\Uc}}

\newcommand{\UU}{{\mathsf{U}}}
\newcommand{\Uc}{\UU_c}

\newcommand{\HH}{\mathsf{H}}
\newcommand{\Hc}{\HH_c}

\newcommand{\cc}{w}

\newcommand{\h}{{\mathfrak{h}}}

\newcommand{\sset}{\subset}
\newcommand{\hreg}{{\mathfrak{h}}^{\operatorname{reg}}}
\newcommand{\ureg}{{\mathsf{U}}^{\operatorname{reg}}}
\newcommand{\Hreg}{{\mathsf{H}}^{\operatorname{reg}}}

\newcommand{\D}{{\scr D}}
\renewcommand{\P}{{\mathbb{P}}}
\newcommand{\X}{{\f X}}

\newcommand{\into}{{\,\hookrightarrow\,}}

\newcommand{\G}{\Gamma}

\newcommand{\Z}{{\mathbb{Z}}}

\newcommand{\Id}{\operatorname{Id}}
\newcommand{\Hilb}{{\operatorname{Hilb}^n{\mathbb{\C}}^2}}

\newcommand{\Coh}{{\operatorname{Coh}}}
\newcommand{\spher}{{\operatorname{spher}}}

\newcommand{\CF}{{\mathcal F}}

\newcommand{\Vo}{{V^\circ}}

\newcommand{\ch}{{\mathsf{ch}}}

\newcommand{\reg}{{\operatorname{reg}}}

\newcommand{\GG}{{\mathfrak{G}}}
\newcommand{\de}{{\D ({\GG})}}
\newcommand{\DD}{{\mathfrak{D}}}
\newcommand{\sgn}{{\sf sign}}
\newcommand{\hr}{{\mathfrak{h}^{\text{reg}}}}

\newcommand{\GGr}{\mathfrak{G}^{\text{reg}}}
\newcommand{\tstar}{\times}
\newcommand{\defn}{\stackrel{\mbox{\tiny{\sf def}}}=}

\begin{document}
\title{Differential operators and Cherednik algebras}
\author{V. Ginzburg}\address{(VG) Department of Mathematics, University of Chicago,  Chicago, IL 
60637, USA.}
\email{ginzburg@math.uchicago.edu}
  \author{I. Gordon} \address{(IG) School of Mathematics and Maxwell Institute for Mathematical 
Sciences,
 Edinburgh  University, Edinburgh EH9 3JZ, Scotland.}
\email{i.gordon@ed.ac.uk}
 \author{J. T. Stafford}
\address{(JTS) Department of Mathematics, 530 Church Street, University of Michigan, Ann Arbor,
MI 48109-1043, USA and  Department of Mathematics, The University of
Manchester, Manchester, M13 9PL,
UK.}
\email{jts@umich.edu,  Toby.Stafford@manchester.ac.uk}
   \thanks{The first    author was  supported in part by the NSF  through the  grant  DMS-0601050.
   The second author was supported by the Leverhulme Trust. The third
author is a Royal Society-Wolfson Research Merit Award holder  and was supported in part by the  NSF  
through   the grant
DMS-0555750 and by the Leverhulme Trust.} \keywords{Cherednik
algebra, Hilbert scheme, characteristic varieties}
\subjclass{14C05, 32S45,16S80, 16D90, 05E10} 

 
\begin{abstract}   
We establish a link between
two geometric approaches to the representation theory of rational
Cherednik algebras of type $\mathbf A$: one based on a
noncommutative Proj construction \cite{GS}; the other
involving quantum hamiltonian reduction of an algebra of differential operators \cite{GG}.

In the present paper, we combine these two points of view 
by showing that  the process of hamiltonian reduction intertwines a naturally
 defined geometric twist functor on $\D$-modules with
 the shift functor for the Cherednik algebra.
That enables us to give a  direct   and relatively short  proof of the key  result
\cite[Theorem~1.4]{GS} without recourse to  Haiman's  deep results on the 
$n!$ theorem \cite{Ha1}. We also show that the characteristic cycles defined independently in these 
two approaches are equal, thereby confirming a conjecture from \cite{GG}.
\end{abstract}

\maketitle
 \tableofcontents

\section{Introduction}  
\label{mmain}  
\subsection{}  Throughout, $W=S_n$ will denote  the $n^{\mathrm{th}}$ symmetric group for some $n\geq 2$.
For a
 parameter  $c\in \C$  we write  $\HH_c$  for \emph{the Cherednik algebra of type $W$}
 with  \emph{spherical subalgebra}  $\Uc= e\HH_c e$, where  
 $e = \frac{1}{|W|}\sum_{w\in W} w\in   \HH_c$ is  
 the trivial idempotent.    (The formal definitions of this and related concepts   can mostly be found in   Section~\ref{notation-section}.) 
  
  Cherednik algebras have been influential
in a wide range of subjects, having been used for instance
to answer questions in real algebraic geometry, integrable systems, combinatorics, and
symplectic quotient singularities.   They are  also closely related to the Hilbert scheme $\Hilb$ of points in the 
plane, a connection  that was formalized in   \cite{GS}  and  \cite{KR}, where it was shown that one can describe  
$\UU_c$ as a kind of  noncommutative deformation of $\Hilb$.
 This was then used in \cite{GS,GS2} to apply  the geometry of $\Hilb$ and 
 Haiman's work on the $n!$ conjecture \cite{Ha1} to analyse the 
 representation theory  of   $\HH_c$.  In the process one finds that  important classes of 
 representations of the Cherednik algebra correspond to important classes of 
 sheaves on the Hilbert scheme and, through them, to important combinatorial objects. 
 
\subsection{}    There are two main aims in this paper. First, 
 the proof of the main result  from \cite{GS}, Theorem~1.4, is  heavily dependent on \cite{Ha1} 
and so one would like to obtain a proof of that result that is independent of  Haiman's work, 
not least because one may then be 
able to apply Cherednik algebra techniques to the study of $\Hilb$ and related combinatorial objects. 
Second, there are two different ways to relate representations of the Cherednik algebra $\HH_c$ to geometric constructions 
and in particular to the Hilbert scheme: 
 through the noncommutative geometry approach of \cite{GS} mentioned above; and through    
  quantum hamiltonian reduction  of the
$GL_n(\C)$-equivariant space $ \gl(\C)\times \C^n$, as in \cite{GG}. One would like to understand the relationship
   between the two approaches. 
  
  In this paper we solve both these problems by using quantum hamiltonian reduction to give an alternative 
  and shorter proof of   \cite[Theorem~1.4]{GS} without
recourse to   \cite{Ha1}.  This also clarifies the relationship between the two approaches, for example showing that  the characteristic cycles defined independently in\cite{GS2} and \cite{GG}  
are actually equal, thereby confirming a conjecture of \cite{GG}.
  
  Before we describe these results in detail  we need to introduce some notation. 

\subsection{} \label{C-defn}  
Let $\h =  \C^n$ denote  the permutation representation of  $W$ and write
  $\hr =\h \sminus \delta^{-1}(0)$
  where  $\delta = \prod_{i< j} (x_i-x_j) \in \C[\h]$ is the discriminant; equivalently $\hr $
   is  the subset of $\h$ on which $W$ acts freely.
 The algebra $\Uc$ will be identified   through  the Dunkl embedding with  a subalgebra of
$\D(\hr)\ast W$, the skew group ring of $W$ with coefficients in the ring
of differential operators on $\hr.$  

Set 
$$\mathcal{C} = \{c\in\C : c=\frac{a}{b}\ \text{where}\  a,b
\in\mathbb{Z} \ \text{with}\ 2\leq b\leq n\}.$$
A scalar $c\in \C$ is called \emph{good} provided that
$c\not\in \mathcal{C}\cap(-1,0)$. The point of the definition is that the good values of $c\in \C$ are the ones  for 
which \cite{GS} and \cite{BE} show that  $\Uc$ has
pleasant properties (see also Remark~\ref{half-integer}).

\subsection{} \label{PQsection}For   $a\in\C$, set  
$$ \PPP{a}{a-1} =  e\mathsf{H}_{a}\delta e \qquad\text{and} \qquad
 \QQQ{a-1}{a} = e\delta^{-1} {\sf H}_a e.$$  By induction, for 
$a\in b+\mathbb{Z}_{\geq 2}$,  define  
 \beq{PQ}
 \PPP{a}{b} = (\PPP{a}{a-1})\cdot(\PPP{a-1}{b})\ \quad \text{and}\ \quad
\QQQ{b}{a} = (\QQQ{b}{b+1})\cdot (\QQQ{b+1}{a}).
\eeq 
   In these equations, the multiplication  is taken inside $
\D(\hreg)\ast W$
 and this makes  both $\PPP{a}{b}$ and $\QQQ{a}{b}$  into   $(\UU_a, \UU_b)$-bimodules.

The key to  \cite{GS,GS2} is the  construction of a  \emph{$\mathbb{Z}$-algebra} $B=\bigoplus_{i\geq j\geq 0}
\bigl(\PPP{c+i}{c+j}\bigr)$ endowed with a natural matrix multiplication. The ring $B$  has a natural 
filtration induced from the differential operator filtration $\Gamma$ on  $\D(\hr)\ast W$ and   
the main result \cite[Theorem~1.4]{GS} showed that {\it if $c+i$ is good for all $i\in \mathbb{N}$,
 then the associated graded ring $\gr_\Gamma B$ of  $B$ is   the $\mathbb{Z}$-algebra associated 
 to $\Hilb$.}  This provides the bridge between Cherednik algebras and Hilbert schemes.

\subsection{} 
There is another way of passing from $\HH_c$ to a more geometric setting which uses the 
hamiltonian reduction from \cite{GG}.  This works as follows.
 Write $V=\C^n$ and set $G=GL(V)$
with Lie algebra $\g=\mathfrak{gl}(V)$.  If  ${\GG} = \g\times V$ with its natural $G$-action, then 
\cite{GG} shows that $\UU_c \cong \left(\D(\GG)/\D(\GG)\cdot I_{c+1}\right)^G$ for an appropriate ideal 
$I_{c+1}$ of $U(\g)$ (see Section~\ref{hamilton-defn} for the details).
One of the main results of this paper shows that there is a natural interpretation of the $\QQQ{a}{b}$ in 
terms of this hamiltonian reduction. Indeed, set $\DD_{c+1} = \D(\V)/\D(\V)\cdot I_{c+1}$ and, 
for $m\in \mathbb{Z}$, consider the  space of semi-invariants
  $$\DD_{c+1}^{\det^{-m}}= \{ D\in \DD_{c+1} : g\cdot D = \det(g)^{-m}D\}.$$

 \subsection{}
  It is easy to check that
 $\DD_{c+1}^{\det^{-m}}$ is a $({\sf U}_{c-m}, \ehe)$-bimodule and our
  first main result shows that it is a familiar object:

\begin{theorem}\label{iso} 
Fix $c\in\C$ and an integer $m\geq1$ 
 so that  each of  $c-1, c-2,\ldots, c-m+1$
 is good {\rm(}this is automatic if $m=1${\rm)}.  Under the differential operator filtration on the two sides
    there is a filtered $({\sf U}_{c-m}, \ehe)$-bimodule  isomorphism 
\begin{equation}\label{iso-equ}
\Theta_{c,m}:  \DD_{c+1}^{\det^{-m}}  
 \iso
\QQQ{c-m}{c}.  
\end{equation} \end{theorem} 
 
 This also leads to a description of  $\PPP{c}{c-m}$ in terms of  the  $ \DD_{d+1}^{\det^{-m}}  $ (see Lemma~\ref{twist2}).
 
 The idea behind the proof of this theorem is   easy to describe. One first shows  that, like $\QQQ{c-m}{c}$, 
 the $({\sf U}_{c-m}, \ehe)$-bimodule  $  \DD_{c+1}^{\det^{-m}}  $ is naturally embedded into 
$\ureg=\Uc[\delta^{-2}]$. Moreover,  both of these bimodules are reflexive on at least one side
(see Corollary~\ref{proj1} and Lemma~\ref{sign-reps}). The theorem is then proved by showing that 
  such a bimodule is unique (see Theorem~\ref{trick1}).
 
 \subsection{} 
Set $A^1= \C[\h\times \h^*]^{\sgn}$, with powers $A^m$ under multiplication in 
 $\C[\h\times \h^*]$. Using Theorem~\ref{iso} we are also able to 
   to give a direct and relatively short  proof of the following result, which is essentially
    \cite[Theorem~1.4]{GS} and is one of the central results from that paper.

\begin{theorem} \label{main-thm}   Keep the hypotheses of Theorem~\ref{iso};
thus    $c-1,  \ldots, c-m+1,$
are all good.
\begin{enumerate}
\item  Under the differential operator filtration $\Gamma$, one has
 $\gr_\Gamma(\QQQ{c-m}{c}) = \delta^{-m}A^me$
inside $\gr_\Gamma\ureg  = \C[\hreg\times \h]e$.

\item  Similarly, 
$\gr_\Gamma (\PPP{c}{c-m}) =  A^m\delta^m e$.
\end{enumerate}\end{theorem}
  
  Crucially,  and unlike the original proof of   \cite[Theorem~1.4]{GS}, 
    the proof of Theorem~\ref{main-thm} does not depend upon Haiman's work \cite{Ha1}.

\subsection{}\label{1.7}    As was remarked earlier,  passing from  
 finitely generated filtered $\HH_c$-modules
to coherent sheaves on $\Hilb$ via  the $\mathbb{Z}$-algebra $B=\bigoplus_{m\geq 0} \PPP{c+m}{c}$
provides a powerful technique for studying the representation theory of $\Uc$ and $\Hc$, see \cite{GS,GS2}. 
Under this relation,  \cite[Corollary ~1.9]{GS} shows that the regular representation ${\sf H}_c$ 
corresponds to  the Procesi bundle $\mathcal{P}$: this is a vector bundle   of rank $n!$ over
 $\Hilb$ whose fibres carry the regular representation of $S_n$, see
 \cite{Ha1}.  These fibres describe Macdonald polynomials and have deep connections with
many areas of representation theory and algebraic combinatorics. 
  
  One may hope to use Theorem~\ref{main-thm}  and the representation
theory of ${\sf H}_c$ to provide a new 
construction of $\mathcal{P}$ which will explain many of its
properties. This goal is still out of reach. Furthermore,
the more detailed  relationship between the category $\Lmod{\HH_c}$ of finitely generated left $\HH_c$-modules  
and the category $\Coh (\Hilb)$ of coherent sheaves on $\Hilb$ established in \cite{GS2} is
 itself dependent upon several important properties of the construction of $\mathcal{P}$ in \cite{Ha1}.
  So  there is much  to understand.

 \subsection{}  In the second part of the paper we use  Theorems~\ref{iso} and \ref{main-thm} to relate other 
 structures appearing in the   $\mathbb{Z}$-algebra and quantum hamiltonian reduction  constructions. 
 These are concerned  with   the {\it functor of hamiltonian reduction}
  $$\widetilde{\BH}_c: \Lmod{(\D_{-c}(\X), SL(V))}\longrightarrow \Lmod{\Uc}$$ 
 where $\X = \g\times \mathbb{P}^{n-1}$ and $\Lmod{(\D_{-c}(\X), SL(V))}$ denotes the  category of
 finitely generated $SL(V)$-equivariant  $\D$-modules on $\X$ that are $nc$-twisted  on $\P^{n-1}$.
The functor $\widetilde{\BH}_c$ is defined formally in \eqref{commute1} but, up to a shift, is given by 
$ \CF\mto \CF^{SL(V)}$.   This  functor has been used  for example in \cite{FG} to relate 
representations of $\Uc$ with generalisations of Lusztig's  character sheaves. 

 Working with 
$\D(\X)$ rather than $\D(\GG)$ is not  particularly significant since one can pass from the latter to the 
former by taking invariants under  the central subgroup $\C^\tstar\subset GL(V)$. 
What is significant  is that there is now   a natural translation    functor 
 $\BS^m : \Lmod{(\D_{-c} (\X), SL(V))} \to  \Lmod{ (\D_{-c-m} (\X), SL(V))}$ 
 given by tensoring with   the  sheaf 
 $\mathcal{O}_{ \mathbb{P}(V) }(nm)$  on $ \mathbb{P}(V)$. On the other hand, 
one has the translation functor $\PPP{c+m}{c}\otimes (-): \Lmod{\Uc} \to \Lmod{\UU_{c+m}}$; 
when $c,c+1,\cdots,c+m$ are good this 
  is an equivalence of categories  that 
   plays a crucial r\^ole in \cite{GS, GS2}, analogous to the
  translation principle in Lie theory.  
 As we prove, it also has a natural 
 description in terms of $\widetilde{\BH}_c$ and its  left adjoint~${}^\top\widetilde{\BH}_c$.

\begin{theorem}  \label{shiftthm} Assume that $n>2$ and that  $c\in\C\smallsetminus \mathbb{Q}_{<0}$.
  Then, for all  integers $m\geq 0$,   there is an isomorphism of functors
  $\widetilde{\BH}_{c+m}\circ \BS^m \circ {}^\top\widetilde{\BH}_c(-)\cong
\PPP{c+m}{c}\otimes_{\UU_{c}} (-)$  that makes
the following diagram commute
$$
\xymatrix{ \strut {\Lmod{(\D_{-c}(\X), SL(V))}} \ 
\ar[rr]^<>(0.5){\BS^m}&& \ 
{\Lmod{(\D_{-(c+m)}(\X), SL(V))}}\ar[d]^<>(0.5){\widetilde{\BH}_{c+m}}\\   
\Lmod{\UU_{c}} \ar[u]^<>(0.5){{}^{\top}\widetilde{\BH}_{c}}   
 \  \ar[rr]^<>(0.5){\PPP{c+m}{c}\otimes_{\UU_{c}}(-)}&& \ \Lmod{\UU_{c+m}}
}
$$ 
\end{theorem}

\subsection{}  A useful tool in the study of Cherednik algebras, just as for Lie algebras, is the concept of 
the characteristic cycle of a $\Uc$-module (this is the same as the characteristic variety, except that it counts
 multiplicities of the irreducible components).
 In fact there are two completely  different constructions of characteristic 
cycles of $\Uc$-modules on $\Hilb$:   the first, $\ch^{GS}$, 
uses  the $\mathbb{Z}$-algebra approach form \cite{GS};  the second, 
$\ch^{GG}$, is defined using the machinery of hamiltonian reduction.  In our final result we prove that
 these two constructions agree, thereby confirming a conjecture from  \cite[(7.17)]{GG}.

\begin{theorem} \label{ccthm} Assume that $n>2$ and that  
 $c\in\C\smallsetminus \mathbb{Q}_{<0}$.
Then for any finitely generated  $\Uc$-module $M$ one has an  equality
of algebraic cycles $\ \ch^{GS}(M)=\ch^{GG}(M)$.
\end{theorem}

\subsection{} The paper is organised as follows. In Section~\ref{notation-section} we introduce the 
basic notation and background material. Section~\ref{hamilton-section}
 is the key to the  paper: it gives a uniqueness result for 
reflexive $(\UU_{c+m},\,\Uc)$-bisubmodules of $\D(\hr)\ast W$.  
We use this in Section~\ref{hamiltonian-section} to   prove Theorem~\ref{iso} in the special case when $m=1$,
 while in Section~\ref{Z-section} we extend this to prove Theorems~\ref{iso} 
and~\ref{main-thm} in general. In Sections~\ref{shift-section} and~\ref{char-section} we prove 
slightly stronger versions of
Theorems~\ref{shiftthm} and~\ref{ccthm}, respectively, which also include the case $n=2$.  
Finally in the appendix we give a detailed 
 proof of the version of hamiltonian reduction that we need, since it does not follow directly from that in \cite{GG}.

\section{Notation and hamiltonian reduction}\label{notation-section}

\subsection{Differential operators}
 \label{diff-subsect}
Let $G$ be a reductive
  algebraic group with Lie algebra  $\mathfrak{g}=Lie(G)$ and write $U(\g)$ for the enveloping algebra
of $\g$.  Let $X$ be  a smooth affine algebraic variety with coordinate ring  $\C[X]$ and 
 ring of regular algebraic differential operators $\D(X)$.  Assume that $G$ acts algebraically on $X$.
 This  gives rise to a  locally-finite $G$-action on $\C[X]$ and  $\D(X)$ via the formul\ae:
\begin{equation}\label{group-action}
(g\cdot  f)(x) \defn f(g^{-1}\cdot x)\quad \text{ and }\quad (g\cdot \theta)(f) \defn  g\cdot (\theta(g^{-1}\cdot f)), 
\end{equation} 
for 
$g\in G, $ $ f\in \C[X],$ $  \theta\in \D(X)$  and $x\in X.$

We should emphasise that this is not the action used in \cite{GG} and \cite{BFG}; those papers implicitly
use the rule  $(g\cdot f)(x) = f(g\cdot x)$ for 
$g\in G, $ $ f\in \C[X]$  and $x\in X.$   (See \cite[Equation~A4]{GG}  and the comments after  
\cite[Lemma~5.3.3]{BFG}.)

 The  action of $G$ on $\C[X]$ and  $\D(X)$  is by given algebra automorphisms, and we let  $\D(X)^G
\sset \D(X)$
denote the subalgebra of $G$-invariant differential operators.
 The action of $G$  on  $X$  gives rise to a $\g$-action on $\C[X]$ by  derivations
 and this induces a Lie algebra map $\tau: \mathfrak{g}\to \mathsf{Der}(\C[X])\subset \D(X)$. 
This   extends uniquely to an associative algebra morphism
$\tau: U(\g)\to \D(X)$.
  
\subsection{}\label{2.2}  Given  a group character $\chi: G\to\C^\tstar$ and 
 a $G$-module $M$,  write  $$M^\chi =\{m\in M\,\mid\, g\cdot m=\chi(g)
m,\en\forall g\in G\}$$ for the corresponding $\chi$-isotypic
component.   Abusing notation, we also write
$\chi:\g\to\C$ for the differential  of the group character $\chi$.

Let $\nu: \g\to\C$ be a Lie algebra character and write  $I_\nu$ for the two-sided ideal in $U(\g)$
generated by the elements $\{x-\nu(x)\ :\ x\in\g\}$.
It is standard that multiplication in $\D(X)$ induces an  algebra structure on
 the space of $G$-invariants $[\D(X)/\D(X)\tau(I_\nu)]^G$, \cite[Section~3.4]{BFG},
called the 
  \emph{quantum hamiltonian reduction of $\D(X)$ at $\nu$}.  Similarly, for a character $\chi$ of $G$, the
   next lemma shows that 
we obtain a natural bimodule structure for the isotypic component  $ \left[\D(X)/\D(X)\tau(I_\nu)\right]^\chi$.
For notational simplicity we will often write $\D(X)I_\nu$ in place of $\D(X)\tau(I_\nu)$.

    \begin{lemma}\label{sublemma1} Let $\chi : G\longrightarrow \C^\tstar$ be a group character and $\nu:\g
\longrightarrow \C$ a Lie algebra character. Multiplication in the algebra
    $\D(X)$ endows     $ \left[\D(X)/\D(X)\tau(I_\nu)\right]^\chi$ 
 with   the structure of a  right $[\D(X)/\D(X)\tau(I_\nu)]^G$-module
   and a left $[\D(X)/\D(X)\tau(I_{\nu+\chi})]^G$-module. 
   \end{lemma}

      \begin{proof} 
The  right  $[\D(X)/\D(X)\tau(I_\nu)]^G$-module   structure  is obvious.
 The differential of the $G$-action on $\D(X)$ induces an action of $\tau(\g)$ by commutation: thus if 
 $x\in \g$ and $D\in \D(X)^\chi$ then $[\tau(x),\,D]=\chi(x)D$. 
 Hence 
    \begin{equation*}\label{righttoleft}
     \bigl(\tau(x)- (\nu+\chi)(x)\bigr)D = D\bigl(\tau(x) -  \nu(x)\bigr)\in \D(X) \tau(I_{\nu})\end{equation*} 
     and so  the left action of
   $\D(X)^G$ on  $ \left[\D(X)/\D(X)\tau(I_\nu)\right]^\chi$  
    factors through the factor ring  $\left[\D(X)/\D(X) \tau(I_{\nu+\chi})\right]^G$.
\end{proof}

\subsection{Rational Cherednik algebras}\label{cherednik-defn}
Fix a positive integer $n\geq 2$, let $W={S}_n$ be  the
symmetric group and write    $e = \frac{1}{|W|}\sum_{w\in W} w\in \C W$ for  the 
trivial idempotent. Similarly, let  $e_-  = \frac{1}{|W|}\sum_{w\in W} \sgn(w)\cdot w$
denote the sign idempotent. Recall that $\h = \C^n$ is the
permutation representation of $W$ and that 
   $$\hr = \{ (z_1, \ldots , z_n)\in\h : z_i\neq z_j \text{ for all }1\leq i < j\leq n\}$$
  denotes  the subset of $\h$ on which $W$ acts freely. Equivalently, if $\{x_1,\dots, x_n\}$ is the basis of
   $\h^*\subset \C[\h]$ consisting of coordinate functions, then   $\hr =\h \sminus \delta^{-1}(0)$
  where  $\delta = \prod_{i< j} (x_i-x_j) \in \C[\h]$ is the discriminant.  There is an induced action of 
  $W$ on $\D(\hreg)$  and we write $\D(\hreg)\ast W$ for  the corresponding skew group ring; for
  $D\in \D(\hreg)$ and $w\in W$ multiplication is
 defined by $wD= (w\cdot D) w$.

Fix a scalar $c\in \C$ and let $\Hc$ denote the \emph{rational Cherednik algebra 
 corresponding to $GL_n(\C)$}. As in \cite[Proposition~4.5]{EG}, 
   we will identify $\Hc$ with  the subalgebra of  $\D(\hreg)\ast W$ generated by $W$, 
the vector space $\h^*=\sum x_i\C\sset\C [\h]$ of linear functions,  and the \emph{Dunkl operators}
\begin{equation} \label{dunkdef}D_c(y_i) \ = \  \frac{\partial}{\partial x_i} \ -\  \frac{1}{2}\sum_{j \neq k} 
c\frac{\langle y_i, x_j-x_k\rangle}{x_j - x_k} (1 - s_{jk})\end{equation} 
where $\{y_1, \ldots , y_n\} \subset \h$ is  the dual basis to $\{x_1, \ldots , x_n\}$
and  the $s_{jk}\in W$ are simple trans\-positions.  By   \cite[Theorem~1.3]{EG} there is a
 PBW isomorphism  $\Hc\cong \C[\h]\otimes \C W\otimes \C[\h^*]$  of $\C$-vector 
spaces. The  \emph{spherical subalgebra} of $\Hc$ is  $\Uc = e\Hc e$.

 Clearly  
$\delta^2\in\C[\h]^W.$ Let $\ureg = \ehe[\delta^{-2}]$ and $\Hreg =
\Hc[\delta^{-2}] =\Hc[\delta^{-1}]$  be the corresponding localised algebras.
By definition,  $\D(\hreg)=\D(\h)[\delta^{-1}]$ and, by  \cite[Proposition~4.5]{EG},
 $\ureg= e\D(\hreg)^We= e\bigl(\D(\hreg)\ast W\bigr)e$ is independent of the choice of~$c$.

  \medskip  \noindent
 \subsection{Remark}\label{half-integer}  (i)
 If $\hccp$ is the corresponding  $SL_n(\C)$
  version of the Cherednik algebra  with spherical subalgebra $\Uc'$ 
  then $\Hc\cong \hccp\otimes_\C \D(\mathbb{A}^1)$ and $\Uc\cong \Uc'\otimes \D(\mathbb{A}^1)$.
It is therefore straightforward to apply the results of (for example) \cite{GS} to $\Hc$ and $\Uc$.
  
 (ii)  The   results in \cite{GS} also assumed that 
$c\notin \frac{1}{2}+\mathbb{Z}$, but this  condition has since been removed by \cite{BE}. So all the
 results in \cite{GS,GS2} can now be applied without that hypothesis. Thus, for example, 
 \cite[Corollary~3.13]{GS} and \cite[Theorem~4.1]{BE} show that if  $c\in \C$ is good,
 then  $\Uc$ is Morita equivalent to $\Hc$ and consequently has finite homological global dimension.

\subsection{}\label{filtration-defn}
The order of differential operators induces a 
filtration on $\D(\hreg)\ast W$ by putting $W$ into degree zero.  
 Essentially every (noncommutative) ring $R$ or module $M$ that we consider  is 
 naturally embedded as a subfactor of either $\D(\hreg)\ast W$ or $\D(X)$ for some variety $X$: we will
  call the induced filtration on $R$ and $M$ 
the \emph{differential operator filtration} 
and this will usually be written as  $M=\bigcup_{j\geq 0}  \Gamma_iM$.
The PBW isomorphism can then be rephrased as saying that there are algebra isomorphisms
 $\gr_\Gamma \Hc\cong \C[\h\times\h^*]\ast W$  and
$\gr_\Gamma \Uc\cong \C[\h\times \h^*]^W$.

\subsection{Quantum hamiltonian reduction}\label{hamilton-defn}
 Fix an $n$-dimensional $\C$-vector space $V$ and put $G=GL(V)$
with Lie algebra $\g=\mathfrak{gl}(V)=\text{Lie}(G)\supset
\mathfrak{sl}(V)$. Write ${\GG} = \g\times V$ with the $G$-action 
 $g \cdot (X,v) = (gXg^{-1}, gv)$ for $g\in G$ and $ (X,v)\in \g\times V$.   
Let $\eu\in\g$ denote the identity. For any $c\in \C,$
let $\chi_c : \g \to \C$ be the Lie algebra homomorphism
 defined by $x\mapsto c\cdot \Tr(x)$. 
For simplicity we write $I_c=I_{\chi_c}\sset U(\g)$: thus
$I_{c}=U(\g)  \mathfrak{sl}_n+U(\g) (\eu-nc)$.
Much of the paper will be concerned with the objects
 \begin{equation}\label{DD-defn}
 \DD_c \ \defn \  \frac{\de}{\de  \tau(I_{c})}  \qquad\text{and}\qquad 
  \DD_{c}^G \ \defn \ \left[\frac{\de}{\de \tau(I_{c})}\right]^{G},
\end{equation}
which, following the earlier convention, are given the differential operator filtration
 $\Gamma$ induced from that on $\D(\GG)$.
 
\subsection{} \label{hamilton-defn2} The action of $\eu$ on $\C[\GG]$ will be  used a number 
of times in this paper and so it is appropriate 
to be explicit about it.  The action of the centre  $\C^\tstar$ of $G$  on
$\GG$ is given by dilation in the second component:
 $\lambda \cdot  (X,v) = (X, \lambda v)$ for $\lambda\in \C^\tstar$ and $(X,v)\in \g\times V =\GG$.
By \eqref{group-action}, $\C^\tstar$ therefore  acts on $V^*\subset \C[V]$ by anti-dilation. 
Since the action of $\g$ is the differential of the $G$-action, this implies that  the action of 
$\tau(\eu)$ will be concentrated purely on $\C[V]$ and that it will then be the negative of the 
Euler operator. In other words, if $\{e_i\} $ is a basis of  $V^*\subset \C[V]$ then 
$\tau(\eu) = - \sum_{i=1}^n e_i\frac{\partial}{\partial e_i}\in \D(V)\subset \D(\GG).$

\subsection{} One of the main results of \cite{GG} shows that 
$\Uc$ may be obtained from  $\D(\GG)$ via  quantum
hamiltonian reduction.

\begin{theorem}\label{first-identity}
With the differential operator filtrations described above, there is for every $c\in \C$ an 
 isomorphism of filtered algebras  $\DD_{c}^G \cong \UU_{c-1}.$  This induces an isomorphism of 
 graded algebras  $\gr \DD_c^G\cong \gr \UU_{c-1}\cong \C[\h\times \h^*]^W$.
\end{theorem}

\begin{proof}   
 This result is a variant of   \cite[Theorem~1.5]{GG}, but since the result does not follow directly from the results
  in that paper, we give a complete proof in the appendix.
  \end{proof}
 

\section{Uniqueness of reflexive bimodules}\label{hamilton-section}

\subsection{} The way we will prove the isomorphism $\DD_{c+1}^{\det^{-m}}  
 \iso \QQQ{c-m}{c}$ of Theorem~\ref{iso} 
 is to note that both sides
are isomorphic to  $(\UU_{c-m},\, \Uc)$-bisubmodules of $\D(\hreg)\ast W$ that are reflexive 
$\UU_{c-m}$-modules. The isomorphism will then follow once we know that such a bimodule is unique.  
The aim of this section is to prove such a uniqueness result, but since the idea  works for more than just 
spherical algebras we will prove it under the following general hypotheses. In what follows we write 
$\gk(M)$ for the Gelfand-Kirillov dimension of a module $M$.

\subsection{Hypotheses}\label{reflexive-hypotheses}
We assume that $(S, \Gamma)$ is a filtered algebra over a  field $k$  such that $\gr_\Gamma S$ is a 
commutative domain. Assume that $R_1$ and $R_2$ are two   subalgebras of $S$ such that:
\begin{enumerate}
\item Each  $R_i$ is a Goldie domain with $S$  contained in the common  Goldie quotient ring $F$ of the $R_i$.
\item Under the induced filtration $\Gamma$, the rings $\gr_\Gamma R_i$
are  Gorenstein algebras that are  finitely generated  modules
 over  a common graded finitely generated $k$-algebra~$C$.
 \item  For each nonzero ideal $I$ of $R_i$ we have $\gk R_i/I\leq \gk R_i-2$.
 \end{enumerate}

\subsection{} We first check that Hypotheses~\ref{reflexive-hypotheses} are satisfied by $\Uc.$

\begin{lemma}\label{jts_prop}  For $i=1,2$,
let  $R_i =a_i\UU_{d_i}a_i^{-1}\subseteq S=\ureg$, for some $d_i\in \C$ and non-zero
$a_i\in \ureg$.  Filter $S$ and its subsets by the differential operator filtration.
Then  Hypotheses~\ref{reflexive-hypotheses} are satisfied by
$R_1$ ad $R_2$.
\end{lemma}

\begin{proof} 
By \cite[Theorem~1.3]{EG} and \cite[Lemma 6.8(1)]{GS},  
$$\gr_\Gamma R_1=\gr_\Gamma R_2 = \gr_\Gamma eH_{d_1}e = 
\C[\h\times \h^*]^W,$$ which is Gorenstein by Watanabe's Theorem \cite{Wa}.
 Thus  parts (1) and (2) of Hypotheses~\ref{reflexive-hypotheses} 
follow, with $C=\C[\h\times \h^*]^W$.

In order to show that part (3)   holds, it suffices to work with $R_1=\Uc$.
Since   $\ureg\cong \D(\h^\reg)^W$,  it  is a simple ring and so 
 $\delta^{2m} \in I$ for some 
 $m>0$. On the other hand, there is a ``Fourier'' automorphism $\kappa$ of $\ehe$
 that maps $\h$ to $\h^*$, see \cite[p.~283]{EG}. Thus,
   the ideal $ \kappa^{-1}(I)$ contains $\delta^{2n}$ for some $n>0$ 
  and so $\kappa(\delta)^{2n}\in I\cap\C[\h^*]$. The PBW isomorphism
 \cite[Theorem~1.3]{EG}  implies that  $\gr \ehe\cong \C[\h\oplus\h^*]^W$
  is  finitely generated  as a module over  its subring $\C[\h]^W\otimes_\C \C[\h^*]^W$. Consequently,
  $\gr (\ehe/I)$ is a finitely generated as a module over the  algebra
  $\C[\h]^W/(\delta^{2m})\otimes_\C\C[\h^*]^W/(\sigma^{2n})$, where
 $\sigma$ is
  the  principal symbol of $\kappa(\delta)$.  This algebra  has
   Gelfand--Kirillov dimension  at most $ 2\dim \h-2$, and hence so does $\ehe/I$.
 \end{proof}

 \subsection{}  Keep the assumptions of Hypotheses~\ref{reflexive-hypotheses}  and set $R=R_1$. 
 Given a non-zero  finitely generated left  $R$-module $M\subset F$ then  
  \cite[Proposition~3.1.15]{MR} shows that there is a canonical identification
$\Hom _{R}(M,\, R)=\{f\in F\mid Mf\sseq R\}$.   The analogous result holds for right modules and the
 \emph{reflexive hull} of $M$ 
is the $R$-module $M^{**}  = \Hom_R(\Hom_R(M,\,R),\,R)\subset F$. Clearly $M\subseteq M^{**}$ 
and $M$ is \emph{reflexive} if $M=M^{**}$.
 Note that when $M=\PPP{c+1}{c}$ or $\QQQ{c}{c+1}$ as modules over either  $R=\Uc$ or $R=\UU_{c+1}$
 these identifications take place inside  $S=\ureg$.

We will need the following application of a theorem of Gabber.
 
 \begin{proposition}\label{gabber} Keep the assumptions of Hypotheses~\ref{reflexive-hypotheses}, set $R=R_1$
and  let $M$ be a  non-zero finitely generated left $R$-submodule of $F$.
 Then $M^{**}$ is the unique largest  left $R$-submodule $M' \subset F$ containing $M$ and such that 
 $\gk M'/M\leq \gk R-2$.
 \end{proposition}
 
 \begin{proof}  Define the {\it grade} of a finitely generated  $R$-module $N$ to be 
 $$j(N) = \min \{ k: \Ext_{R}^k (N, R) \neq 0 \}.$$ By Hypotheses~\ref{reflexive-hypotheses}(2) we may apply
  \cite[Theorems~4.1 and 4.3]{Bj} to conclude that 
$R$ is  AS-Gorenstein and  that $j(N) = \gk R  - \gk N$  for all finitely generated $R$-modules 
$N$. Thus if $M\subseteq M'$ 
are finitely generated $R$-submodules of $F$, then $\gk 
M'/M \leq \gk R-2$ if and only if $j(M'/M) \geq 2$.  Such modules $M'$ are called {\it tame pure 
extensions} of $M$. By \cite[Theorem~3.6]{BjE},  if $M'$ is any tame pure extension of $M$ then there 
exists an injection $\alpha: M' \hookrightarrow M^{**}$ that is the identity on $M$. But since $M'/M$ is a 
torsion $R$-module, for any $m'\in M'$  there exists  non-zero $r\in R$ such that 
$r m' \subseteq M$ and so $r \alpha(m') = \alpha( rm') = rm'$. This proves that $\alpha(m') = m'$ for 
all $m'\in M'$ and  finishes the proof of the proposition.
 \end{proof}

\subsection{}  We are now ready to prove our uniqueness result for reflexive bimodules. 

\begin{theorem}\label{trick1} Assume that  $(R_1,R_2,S)$ satisfy Hypotheses~\ref{reflexive-hypotheses} and
 let $M$ be a non-zero $(R_1,\,R_2)$-bisubmodule of $S $ that is finitely generated 
and reflexive on one side. Then it is reflexive and finitely generated on the other side and is the unique 
such object.
\end{theorem}

\begin{proof} By symmetry, we may assume that $M$ is a finitely generated reflexive left $R_1$-module;
thus  $M=M^{**} \defn \Hom_{R_1}(\Hom_{R_1}(M,\,R_1),\, R_1)$. 
Part (1) of Hypotheses~\ref{reflexive-hypotheses}  implies that, as a left $R_1$-module,
 $M\subseteq R_1f$ for some $f\in S $ from which 
we conclude that   $\gr_\Gamma M \subseteq (\gr_\Gamma R_1)\sigma(f)$,  where 
$\sigma(f)$ is the principal symbol of $f$. By Hypotheses~\ref{reflexive-hypotheses}(2)  
  $\gr_\Gamma M$ is therefore a finitely generated
 (left) module over  both $\gr R_1$ and $C$.
But  $\Gamma$ is also a filtration of $M$ as a  right $R_2$-module. Since $\gr_\Gamma M$ is 
a finitely generated right $C$-module it is also finitely generated as a right 
$\gr_\Gamma R_2$-module. Thus $M$ is a finitely generated right $R_2$-module, see \cite[Proposition~6.4]{KL}.
We remark that it  now  follows from \cite[Corollary~3.4]{KL} that the Gelfand--Kirillov dimension of $M$ 
as a left   $R_1$-module equals
that of $M$ as a right $R_2$-module and we need not distinguish between them.

Next, suppose that $M$ is not unique and that $N$ is a second such bimodule.
 Let $\cmpl{M}=(M+N)/N$; by symmetry we may  assume that $\cmpl{M}\not=0$.
Then $\cmpl{M}$ is a finitely generated left $R_1$-module, say $\cmpl{M}=\sum_{j1}^r R_1
\cmpl{m}_i$.
Since $R_2 $ is an Ore domain, $\cmpl{M}$ is a torsion right $R_2$-module and so
 $$I=\text{r-ann}_{R_2} \cmpl{M} = \bigcap_{j=1}^r \text{r-ann}_{R_2} \cmpl{m}_j
\not=0.$$  By Hypotheses~\ref{reflexive-hypotheses}(3) this implies that $\gk \cmpl{M}\leq\gk R_2/I \leq  \gk R_2 -2.$
Thus we also have $\gk (M+N)/M\leq \gk R_1-2$ as   left $R_1$-modules.  Our hypotheses mean that   
  Proposition~\ref{gabber} can be applied, from which it follows that 
$N \subseteq (M+N)\subseteq M^{**}=M$.
By symmetry,  $M=N$ and so $M$ is indeed unique.

It remains to prove that $M$ is reflexive as a right $R_2$-module, so  suppose to the contrary  that 
$N=\Hom_{R_2}(\Hom_{R_2}(M, \, R_2) , \,  R_2)\supsetneq M$. By Proposition~\ref{gabber}, 
again, it follows that  $\gk N/M\leq \gk R_2-2$.
But as $M$ is a left $R_1$-module, so is $N$ and so by  applying 
Proposition~\ref{gabber} to the left $R_1$-modules $M\subseteq N$ we conclude that $N\subseteq 
M^{**}=M$.
\end{proof}

\subsection{}\label{P-chat} 
In applications of Theorem~\ref{trick1} it is important to be careful to ensure that the bimodule 
structure of $M$ is  induced from that of $S$; after all, for any nonzero $s\in S$ the vector space $R_1s$ 
is an $(R_1,\,s^{-1}R_1s)$-bimodule   and hence, up to isomorphism, an $(R_1,\, R_1)$-bimodule.
  
We are going to apply Theorem~\ref{trick1}  to the modules $ \PPP{c+1}{c}$  and $\QQQ{c}{c+1} 
$ from  Section~\ref{PQsection} and so it is worth emphasizing that  no such problems occur here. Indeed, from
 \cite[Proposition~4.1]{BEG2}, $\Uc=e\delta^{-1}\HH_{c+1}\delta e$ as subrings of $\ureg$. Therefore, 
 the $(\UU_{c+1},\,\Uc)$-bimodule structure of $\PPP{c+1}{c}= e\HH_{c+1}\delta e$ is indeed induced
  from multiplication in  $S=\ureg$.   Similar comments apply to $\QQQ{c}{c+1}$. 
  
\subsection{}   We are now ready apply  
Theorem~\ref{trick1}   to  $ \PPP{c+1}{c}$  and $\QQQ{c}{c+1}$.
\begin{corollary}\label{proj1}
\begin{enumerate}
\item  As a right $\Uc$-module, $\PPP{c+1}{c}$ is projective whenever $c$ is good.
As a left  $\UU_{c+1}$-module, $\PPP{c+1}{c}$ is projective whenever $c+1$ is good.
For all values of $c$, the module $\PPP{c+1}{c}$ is reflexive on both sides.

\item  As a left $\Uc$-module, $\QQQ{c}{c+1}$  is projective whenever $c$ is good.
As a right  $\UU_{c+1}$-module, $\QQQ{c}{c+1}$ is projective whenever $c+1$ is good.
For all values of $c$, the module $\QQQ{c}{c+1}$ is reflexive on both sides.

\item  For all values of $c$,  and under their natural embedding into $\ureg$, we have 
\begin{align*}\PPP{c+1}{c}\ &= \  \Hom_{\Uc}(\QQQ{c}{c+1},\, \Uc)
 \ = \  \Hom_{\UU_{c+1}}(\QQQ{c}{c+1},\, \UU_{c+1}),\quad\text{resp.}\\ 
 \QQQ{c}{c+1}\ &= \  \Hom_{\Uc}(\PPP{c+1}{c},\, \Uc)
\ = \  \Hom_{\UU_{c+1}}(\PPP{c+1}{c},\, \UU_{c+1}).
\end{align*}
\end{enumerate}
\end{corollary}

\begin{proof}  (1,2)  A slightly weaker version of this result is given in  \cite{GS}, but we give a 
different proof since we will need the  argument later.
We will just prove the result for $Q=\QQQ{c}{c+1}$; the same argument works for part~(1).

 Under the differential operator filtration $\Gamma$, the proof of  \cite[Lemma~6.9(2)]{GS}  shows  that 
$ \gr_\Gamma Q   \cong A^1= \C[\h\times \h^*]^{\sgn}$ as modules
over $A^0=\C[\h\times\h^*]^W$.    By the Hochster--Eagon Theorem,
  $\C[\h\times \h^*]$ is a Cohen--Macaulay $\C[\h\times \h^*]^W$-module
  and hence so is its summand $\C[\h\times\h^*]^\sgn$. By
  \cite[Corollary~3.12]{Bj} we deduce that $Q$
    is a torsion-free, Cohen--Macaulay $\UU_{c+1}$-module
   in the sense that $\text{Ext}^{j}_{\UU_{c+1}}(Q,\,  \UU_{c+1})=0$ for $j>0$.
   But, if  $c+1 $ is good,  then $\UU_{c+1}$ has finite global dimension by \cite[Corollary~3.15]{GS} and
    \cite[Corollary~4.3]{BE}  and so this implies that 
    $Q$ is  a projective right $\UU_{c+1}$-module. 
    
   The analogous argument shows that $Q$ is a projective left $\Uc$-module whenever $c$ is good
   and so  for {\it any}  value of $c$ this implies that on at
    least one side $Q=\QQQ{c}{c+1}$ is projective and hence reflexive.  By Theorem~\ref{trick1}
    and Lemma~\ref{jts_prop},  
    $Q$ is then reflexive on the other side.
 
(3) Since $\QQQ{c}{c+1}$ is a finitely generated $(\Uc,\,\UU_{c+1})$-bisubmodule of $\ureg$ that is 
reflexive on both sides, the dual module  
 $\QQQ{c}{c+1}^* =  \Hom_{\Uc}(\QQQ{c}{c+1},\, \Uc)$  is a   nonzero  
 $(\UU_{c+1},\,\Uc)$-bisubmodule  of $\ureg$
 which is reflexive and finitely generated  as a  right $\Uc$-module.   But, by part~(1), the same is true of $
\PPP{c+1}{c}$.   Hence  $\QQQ{c}{c+1}^*  = \PPP{c+1}{c} $ by Theorem~\ref{trick1}     and Lemma~\ref{jts_prop}. 
The  same argument can be used  to  prove the rest of part~(3).
\end{proof}

\subsection{} The following easy and well known result will be used frequently.

\begin{lemma}\label{proj-tensor}  Let $R,S$ be rings with a projective right $R$-module $M$ and
an $(R,S)$-bimodule $N$ that is projective as a right $S$-module. Then 
$M\otimes_RN$ is a projective right $S$-module.
\qed\end{lemma}

\subsection{} It is now easy to generalise Corollary~\ref{proj1}  to the modules $\PPP{c+m}{c}$  and 
$\QQQ{c}{c+m}$ of \eqref{PQ}.

\begin{theorem} \label{the trick}  Fix $c\in\C$ and an integer
$m\geq1$ such  that   $c+1,c+2,\dots, c+m-1$ are all  good.

 \begin{enumerate}
\item 
$\PPP{c+m}{c}$ is the unique   non-zero  $(\UU_{c+m},\, \Uc)$-bisubmodule of $\ureg$ that is reflexive 
as  either  a  right   $\Uc$-module or a left $\UU_{c+m}$-module. Multiplication in $\ureg$  induces an 
isomorphism  of $ (\UU_{c+m},\, \Uc)$-bimodules,
\begin{equation}\label{trick-equ1}
 \ (\PPP{c+m}{c+m-1})\otimes_{\UU_{c+m-1}}\cdots\otimes_{\UU_{c+1}} (\PPP{c+1}{c}) \  \stackrel{\sim} 
\longrightarrow \  \PPP{c+m}{c }.
\end{equation}
 
\item 
  $\QQQ{c}{c+m}$ is the unique nonzero $(\Uc,\, \UU_{c+m})$-bisubmodule of $\ureg$ that is reflexive as 
either  a  right  $\UU_{c+m}$-module or a left $\Uc$-module. Multiplication gives an isomorphism of 
$(\Uc,\, \UU_{c+m})$-bimodules, 
\begin{equation}\label{trick-equ}
 \ (\QQQ{c}{c+1})\otimes_{\UU_{c+1}}\cdots\otimes_{\UU_{c+m-1}} (\QQQ{c+m-1}{c+m}) \ \stackrel{\sim}
\longrightarrow   \ \QQQ{c}{c+m }.
  \end{equation}

\item Either $c$ or $c+m$ is good. In the former case $\PPP{c+m}{c}$ and $\QQQ{c}{c+m}$ are 
projective $\Uc$-modules, while in the latter case they are projective $\UU_{c+m}$-modules.
\end{enumerate}
\end{theorem}

\begin{proof} We only prove assertions for  $\QQQ{c}{c+m}$; the proofs for $\PPP{c+m}{c}$
are essentially the same.  We will also assume that $c+m$ is good; the proof when $c$ is good is again 
very similar.

The hypotheses on $c$ now ensure that, by  Corollary~\ref{proj1},  $\QQQ{c+i}{c+i+1}$ is projective 
 as a right $\UU_{c+i+1}$-module for all $0\leq i\leq m-1$. 
 It follows from Lemma~\ref{proj-tensor} and  a routine induction that \eqref{trick-equ} holds and hence  
that  $\QQQ{c}{c+m}$ is both projective and finitely generated as a
   right $\UU_{c+m}$-module.   The remaining assertions   follow from Theorem~\ref{trick1} and Lemma~\ref{jts_prop}.
\end{proof}

\subsection{Example} \label{Bad-c example}
We end the section by noting that Theorem~\ref{the trick}  does not extend to all
values of $c$. To see this consider the special case when $n=2$ and set
$\UU=\UU_{\frac{1}{2}}$. 
Further, let $\cmpl{U}=U(\mathfrak{sl}_2)/(\Omega)$, where $\Omega$
is the quadratic Casimir element. Using \cite[Example~6.12]{GS2} and 
Remark~\ref{half-integer}, one finds that $\UU\cong \cmpl{U}\otimes_\C \D(\mathbb{A}^1) $.
In this case, $\UU_{-\frac{1}{2}}$ is a simple ring, necessarily equal to $\End_{\UU}({\sf Q})$ for 
${\sf Q}= \QQQ{-\frac{1}{2}}{\frac{1}{2}}$.  Write 
 $$\mathrm{Tr}({\sf Q}) \ = \  \{\theta(q) : \theta \in \Hom_\UU(Q,\UU)\ \text{and}\ q\in {\sf Q}\} \ \subseteq\  \UU$$ for 
the trace
ideal of ${\sf Q}$.    
Since  $\UU_{-\frac{1}{2}}$ is simple but $\UU$ is not,  the projective $\UU$-module 
${\sf Q}$ cannot be a progenerator. Equivalently, ${\sf Q}$ is not projective as a left $\UU_{-\frac{1}{2}}$-module,
 thereby showing that the final sentence of Theorem~\ref{the trick}(3) does not hold for arbitrary   $c$.
 
 Write  $F$ for the field of fractions of  $\UU$  and for a $\UU$-submodule $M\subset F$, identify
 $M^*=\Hom_\UU( M,\,\UU)$ with 
$\{\theta\in F : \theta M\subseteq \UU\}$.
We claim that 
$V=\QQQ{-\frac{3}{2}}{\frac{1}{2}} =(\QQQ{-\frac{3}{2}}{-\frac{1}{2}}) {\sf Q}  $ is not reflexive on either side.
Indeed, suppose that it is reflexive on one side (and hence on both sides by Theorem~\ref{trick1}).
By \cite[Theorem~B]{St}, $\UU\cong \UU_{-3/2}$ has global dimension $2$ and so $V$
 is also projective on both sides.  As $U$ is a maximal order, see \cite[Lemma~3.1]{St}, $\UU_{-3/2}
=\End_{\UU}(V)$ and hence  $V^*=\Hom_{\UU}(V,\UU)$ is also equal to 
$\Hom_{\UU_{-{3/2}}}(V,\, \UU_{-{3/2}})$.
Thus the Dual Basis Lemma applied to $V$ as a $\UU_{-{3/2}}$-module implies that $V^*V=\UU$.
But now $W=V^*(\QQQ{-\frac{3}{2}}{\frac{1}{2}})$ satisfies $W{\sf Q}=\UU$, contradicting the fact that $\sf Q$
 is not a progenerator.  This proves the claim.


\section{Relating hamiltonian reduction and shift functors}\label{hamiltonian-section}

\subsection{}  Recall from \eqref{hamilton-defn} that
$\GG=\g \times V$ 
with the natural action of $G=GL(V)$ and that,  as in \eqref{DD-defn},  we write
$$\DD_{c+1} = \frac{\displaystyle \de}{\displaystyle\de \tau(I_{c+1})};$$
thus  $  \DD_{c+1}^G\cong \UU_{c}$ by Theorem~\ref{first-identity}.
The aim of this section will be to construct  a morphism $\Psi_{c}:\DD_{c+1}^{\det^{-1}}\, \to\, \QQQ{c-1}{c}$
and use it to prove Theorem~\ref{iso} in the case $m=1$. The complete proofs of 
 Theorems~\ref{iso} and~\ref{main-thm}, which are given in Section~\ref{Z-section},
  will then follow by applying results from \cite{GG}.
    We begin with the construction of $\Psi_c$, for which we need to expand upon the 
isomorphism in Theorem~\ref{first-identity}.

   \subsection{}\label{sub}
Let ${\sf vol} \in \wedge^n V^*$ be a non-zero volume element on $V$ and,
  following \cite[Equation~5.3.2]{BFG},  define a map  $\svol: {\GG} \to \C$ by
  \begin{equation} \label{sdefn} \svol(X, v) \  \defn\  \langle {\sf vol}, v \wedge Xv \wedge X^2v \wedge
   \cdots \wedge X^{n-1}v\rangle.\end{equation}
  Clearly $\svol\not=0$ and, as $(g\cdot\svol )(X,v)= \svol(g^{-1}\cdot(X,v)) = \svol(g^{-1}Xg,g^{-1}v)$,
  it is a $\det^{-1}$ semi-invariant under the $GL(V)$-action.\footnote{As we  
  mentioned in  Section~\ref{diff-subsect}, \cite{BFG} and   \cite{GG} use the opposite convention for 
  group actions and so in those papers $\svol$ transforms by the determinant; see, for example, the 
sentence after   \cite[(6.18)]{GG}.} 
Note also that $\svol _{\bigl| \hreg} = \delta$.

  We now want to use the radial component map, described as follows.
  Take the Zariski open dense subset  
    $$\GG\cyc \  \defn\  \bigl\{(X,v):\  v \text{ is a cyclic vector for the  operator $X: V\to V$}\bigr\}\ \subset\  \GG$$
      and, for $c\in \mathbb{C}$, let
  $$\mathcal{O}(\GG\cyc,c\  ) \  \defn \ \{f\in \C[\GG\cyc] :  x\cdot f = c\Tr (x) f \ \  \text{for all} \  x\in \g\}.$$
  By \cite[Corollary~5.3.4]{BFG} restriction of functions  induces an isomorphism
   $\mathcal{O}(\GG\cyc,c) \cong \svol^{-c}\C[\h/W]$.   Moreover,  if $D\in   \D(\GG\cyc)^G$
 then, by restriction, $D$ induces  a differential operator from 
 $\mathcal{O}(\GG\cyc,c)$ to itself. This defines the \emph{twisted radial components}
map 
 \begin{equation}\label{Rad-defn}   
 \Rad_c :\  \D(\GG\cyc)^G\too\D(\h/W),
\qquad D\mto \Rad_c(D) = \svol^{c}\ccirc  \bigl(D|_{\mathcal{O}(\GG\cyc,c)}\bigr)\ccirc \svol^{-c} .
\end{equation}
As in  \cite[Section~6.5]{GG} and in the notation of \eqref{DD-defn},  
  this map is an algebra homomorphism that induces the    isomorphism 
 \begin{equation}\label{sub11}
\Rd_c :\ \DD_c^G =  \left(\D({\GG})/\D({\GG}) \tau( I_{c})\right)^G
 \xrightarrow{\sim} {\UU}_{c-1}\end{equation}
 from Theorem~\ref{first-identity}. The details are given in the appendix to this paper. 
  In particular,   $\Rad_c$ restricts to give a surjection $ \D({\GG})^G \twoheadrightarrow \UU_{c-1}$.

 As is observed after \cite[(6.18)]{GG}, $\GG\cyc=\GG\smallsetminus \svol^{-1}(0)$ and so, if 
  $D\in \D(\GG)^{\det^{-1}}$, then $D\svol^{-1}$ is a well-defined $G$-invariant differential 
  operator on $\C[\GG\cyc]$. Hence  we obtain a map 
   $$\widehat{\Psi}_{c} : \ \de^{\det^{-1}} \longrightarrow\   e \D(\h/W)e \  \subset \  \ureg \ =
    \   e\D(\hr)^We,   \qquad  D\mto  e \Rad(D\svol^{-1}) e .$$
 As in the proof of Lemma~\ref{sublemma1},  we have  $[\tau(x),\svol^{-1}] = \Tr(x)\svol^{-1}$
  for any $x\in \g$ and so  
  $$\bigl(\tau(x) - (c+1)\Tr(x)\bigr)\svol^{-1} =  \svol^{-1}\bigl( \tau(x) - c\cdot\Tr(x)\bigr)
   \in \D(\GG^{\cyc})\tau(I_c).$$ 
  Therefore the map  $\widehat{\Psi}_{c}$ factors through $\DD_{c+1}^{\det^{-1}}$  to give  the desired map
    \begin{equation}\label{sub12}
    \Psi_{c} :\ \DD_{c+1}^{\det^{-1}}\to \ureg.
    \end{equation}
   
 \begin{lemma}\label{sub1}  Under multiplication in $\ureg$,
the image $\mathrm{Im}(\Psi_c)$ is a  non-zero $(U_{c-1},\, U_{c})$-bimodule. 
Moreover,  $\Psi_c$ is a $(\UU_{c-1},\,\Uc)$-bimodule homomorphism in the sense that
 \begin{equation}\label{sublemma2}
\Psi_{c}(\alpha_{c}\beta \alpha_{c+1} ) =
\Rd_{c}(\alpha_{c})\Psi_{c}(\beta)\Rd_{c+1}(\alpha_{c+1})  
\end{equation}
for $\alpha_c\in \DD^G_c$, $ \alpha_{c+1}\in \DD^G_{c+1}$ and $\beta\in \DD_{c+1}^{\det^{-1}}$.
\end{lemma}

\begin{proof}   Note that $\svol \in \D(\GG)^{\det^{-1}}$ and that $\Psi_c([\svol])=e\Rad(\svol^{-1}\svol) e = 
e$. Hence $\Psi_c\not=0$.
Combined with \eqref{sub11}, it now suffices to prove \eqref{sublemma2}.

By Lemma~\ref{sublemma1},  $\alpha_{c}\beta \alpha_{c+1} \in\DD_{c+1}^{\det^{-1}}$ and so the 
 left  hand side of \eqref{sublemma2} is well defined.
So \eqref{sublemma2}     therefore follows from the   computation 
\begin{eqnarray*}
\Rd_{c}(\alpha_{c})\Psi_{c}(\beta)\Rd_{c+1}(\alpha_{c+1})
  & =& \svol^{c} \alpha_{c}\svol^{-c}\cdot \svol^{c}(\beta \svol^{-1})
  \svol^{-c}\cdot \svol^{c+1} \alpha_{c+1} \svol^{-(c+1)} \\
    \noalign{\vskip 5pt}
   &=& \svol^{c} (\alpha_{c}\beta \alpha_{c+1}\svol^{-1})\svol^{-c}
   \ = \ \Psi_{c}(\alpha_{c}\beta\alpha_{c+1}
   ). \end{eqnarray*}
 \vspace{-1.25cm}
 
  \end{proof}

\subsection{}\label{moment}    Let $\mu_{\GG} : T^*{\GG}\to \g^*\cong \g$ be the moment map as
defined in \cite[(2.4)]{GG} and, as in \cite[(1.1)]{GG}, set 
\begin{equation}\label{M}
 \M \  \defn\  \mu_{\GG}^{-1}(0) \ = \  \bigl\{ (X,Y,v,w)\in \g\times\g\times
V\times V^* \mid [X,Y]+vw=0\bigr\}. 
\end{equation}

 Recall that $A^1=\C[\h\times\h^*]^\sgn$. The powers $A^r$ are obtained by 
multiplication inside $\C[\h\times\h^*]$: they will be regarded as  modules over
 $A^0=\C[\h\times\h^*]^W$ under the natural induced structure.
Then, once one recalls our conventions about group actions from
\ref{diff-subsect},  it follows from  \cite[Proposition~A2]{GG}  that 
\begin{equation}\label{A}
  A^r \cong \C[\mathcal{M}]^{\det^{-r}},
\quad r=0,1,2,\ldots.
\end{equation}
This result is stated in \cite{GG} as an isomorphism of 
 vector spaces, but the proof shows that it is in fact an isomorphism of
 graded $A^0$-modules.

\subsection{}  As before,  for any subfactor $T$ of $\de$ with  
the induced differential operator filtration, we write 
$\gr T= \gr_\Gamma T$ for the associated graded object. We will  consider $\C[\h\times\h^*]$
 and its submodules as graded in the second term; equivalently they are given the gradation induced 
 from the identity $\C[\h\times \h^*]=\gr \D(\h)$. Recall that, by Theorem~\ref{first-identity} and  \eqref{sub11}, 
 the   map $\mathfrak{R}_{c}$ induces a graded isomorphism
 $\gr\Rd :\  \gr  \DD_{c}^G \iso  A^0 =  \C[\h\times \h^*]^W.$
 
\begin{lemma}\label{sign-reps} 
{\rm (1)} For all $c\in \mathbb{C}$,  and $r\in\mathbb{N}$,  we have $\gr \DD_{c+1}^{\det^{-r}} \cong A^r$
as graded $A^0$-modules, where   the $A^0$-module structure of $\gr \DD_{c+1}^{\det^{-r}} $ is 
defined  via $\gr \Rd$.

{\rm (2)}   For all $c\in \mathbb{C}$,  $\DD_{c+1}^{\det^{-1}}$ is reflexive as both a right $\Uc$-module 
and a left $\UU_{c-1}$-module.
\end{lemma}

\begin{proof}  (1)
 We note for future reference that the action of $G$ preserves the differential
operator filtration and so, for any $r\in \mathbb{Z}$, we have 
$\gr(\de^{\det^r})=\left(\gr \, \de\right)^{\det^r}$. In particular  we can  identify 
$$\bigl(\gr \DD_{c+1}\bigr)^{\det^{-r}}   \ = \ 
\frac{\gr(\de^{\det^{-r}})}{\gr(\left[\de \tau( I_{c+1})\right]^{\det^{-r}})}$$ without causing ambiguity.

We return to the proof and write $\mu=\mu_\GG$.  By construction there is  a
sequence of graded $\C[\h\times\h^*]^W$-module maps:
\begin{equation}\label{sign-equ}
\begin{array}{rl}
\chi :\  A^r \cong \C[\mathcal{M}]^{\det^{-r}} \ = &
 \left[\frac{\displaystyle \C[T^*{\GG}] }{
\strut \displaystyle\C[T^*{\GG}]\cd \mu^*(\g)} \right]^{\det^{-r}}\ = \
\left[ \frac{\displaystyle\text{gr}\D({\GG}) }{\strut\displaystyle\text{gr}\D({\GG})\cdot
\text{gr}\, \tau(I_{c+1})} \right]^{\det^{-r}} \\ 
\noalign{\vskip 10pt}  
& \ \buildrel{\alpha}\over\ontoo \ 
 \left[ \frac{\displaystyle\text{gr}  \D({\GG}) }{
\strut\displaystyle \text{gr}\bigl(\D({\GG})\tau( I_{c+1})\bigr)}  \right]^{\det^{-r}}
 \end{array}\end{equation}
 By   \cite[Corollary~2.6]{GG}  $\mu$ is flat and so, 
 by \cite[Proposition~2.4]{H},  the final surjection $\alpha$ in \eqref{sign-equ}
   is an  isomorphism.
 Hence $\chi$ is an isomorphism.
   
 (2)  The fact that $\DD_{c+1}^{\det^{-1}}$ is a $(\UU_{c-1},\,\Uc)$-bimodule is part of
  Lemma~\ref{sub1}.  By part~(1), $\gr \DD_{c+1}^{\det^{-1}}  \cong A^1$ as both 
  a right module over $\gr \Uc \cong \mathbb{C}[\h\times\h^*]^W$ and as a left module over
  $ \gr  \UU_{c-1} \cong \mathbb{C}[\h\times\h^*]^W$.  
  The proof of Corollary~\ref{proj1}(2) can now be used unchanged to give the desired result.
\end{proof}

\subsection{}  
We are now able to prove Theorem~\ref{iso} in the case $m=1$ and  we are able to do so
 without restriction on $c$.  Thus we prove:
 
 \begin{proposition}\label{isom=1}
 For any $c\in \mathbb{C}$ the homomorphism $\Psi_c$ from \eqref{sub12} induces an isomorphism 
$\DD_{c+1}^{\det^{-1}}\cong \QQQ{c-1}{c}$ 
of $(\UU_{c-1},\,\Uc)$-bimodules.
\end{proposition}

\begin{proof}  
A  key observation here is that $\Uc$ is a noetherian domain with quotient division ring $F$ containing 
 $\ureg$. Hence any non-zero finitely generated $\Uc$-submodule $M$ of $F$ must be a torsion-free module of 
  Goldie rank one; equivalently,\ every proper factor of $M$  is a torsion module.  

We first claim that   $\Psi_{c}$  is injective.  
 To see this note that, by Lemma~\ref{sign-reps}(1),
  $\gr \DD_{c+1}^{\det^{-1}}$ is a
   torsion-free, rank one  module over the domain  $\text{gr}(\ehe) \cong \C[\h\times \h^*]^W$.
Consequently,   $\DD_{c+1}^{\det^{-1}}$ is  a torsion-free right   $\Uc$-module of Goldie rank one
     and so any proper factor of this module would be torsion (or zero).
   But, by Lemma~\ref{sub1},    $\text{Im}(\Psi_{c})\subset \ureg$
is a non-zero torsion-free $\ehe$-module. Hence
   $\Psi_{c}$  is indeed injective.
 
Therefore, by Lemma~\ref{sign-reps}(2), respectively Corollary~\ref{proj1}(2),
\emph{both}   $\mathrm{Im}(\Psi_c)$ and $\QQQ{c-1}{c}$ are reflexive $(\UU_{c-1},\, \Uc)$-bimodules
of $\ureg$, and we emphasise that in both cases the bimodule structure is that induced from the 
bimodule structure of $\ureg$.
By Theorem~\ref{trick1} and Lemma~\ref{jts_prop} they are therefore equal.
\end{proof}

\subsection{}
We end the section by strengthening the isomorphism from
Proposition~\ref{isom=1} to an isomorphism of filtered spaces.
The difficulty here is that, although 
$\DD_{c+1}^{\det^{-1}}$ and $\QQQ{c-1}{c} $ are given the differential operator filtrations $\Gamma$, 
these are induced  from  two
different  rings of differential operators $\D(\GG)$, respectively $\D(\hreg)$. In the abstract it is not clear
 that the two resulting filtrations on  $\QQQ{c-1}{c}$   are  closely  related (see \cite[Probl\`eme~3.5]{Le}, 
for example).
 Fortunately, this problem is easily resolved in the present setting.

\begin{corollary}\label{iso2} Keep the notation of Proposition~\ref{isom=1}.
Then the isomorphism
$\Psi_{c}$ is a filtered isomorphism in the sense that
$\Psi_{c}(\Gamma_{r} \DD_{c+1}^{\det^{-1}})=\Gamma_{r}(\QQQ{c-1}{c})$ for all $r$.
\end{corollary}

\begin{proof}  Set $D=\DD_{c+1}^{\det^{-1}}$ and $Q=\QQQ{c-1}{c}$.
 The analogous isomorphism of filtered algebras
$\DD_{c+1}^G \cong \Uc$ is given by  Theorem~\ref{first-identity}.
Since the radial component map  $\Rad$   at worst
decreases the degree of an operator, the map $\Psi_{c}$ is at
least a map of filtered modules in the sense that
$\Psi_c(\Gamma_mD) \subseteq \Gamma_mQ$ for all $m \geq 0$. Suppose  that $\Psi_c$ is
not an isomorphism of filtered objects. Since $\Psi_c$ is an isomorphism of unfiltered modules, 
 there must then exist some $ d\in D $ for which the order of $\Psi_c(d)$ is less than the order of $d$.
 Hence  $ \text{gr}(\Psi_c) : \text{gr}_\Gamma D \to \text{gr}_\Gamma Q$ will have a
nonzero kernel. But we know from Lemma~\ref{sign-reps}(1)
 that  $ \text{gr}_\Gamma D$ is a rank one torsion-free $\C[\h\times \h^*]^W$-module, and so  its image
 must therefore be torsion. Since $\text{gr}_\Gamma Q$ is a submodule of the domain
 $\text{gr}_\Gamma \ureg=\C[\hreg\times \h^*]^W$ this forces $\text{gr}(\Psi)=0$.
 This contradicts the fact that, for example,
 $0\not=e = \Psi_c(\svol)  \in \Psi(\Lambda_0D) \subseteq \Gamma_0Q.$
 \end{proof}

\section{Proof of Theorems~\ref{iso} and \ref{main-thm}} \label{Z-section}

\subsection{} 
We continue to write  $\DD_c=\de/\de \tau( I_{c})$  in the notation of \eqref{DD-defn},
  and recall from Theorem~\ref{first-identity} that $ \DD_c^{G} \cong \UU_{c-1}$.
For any $p,q\in \mathbb{Z}_{\geq 0}$, multiplication on $\de$ induces a map
$$\de^{\det^{-p}} \otimes_{\D(\GG)^G} \de^{\det^{-q}} \to \de^{\det^{-(p+q)}}.$$
 It follows from Lemma~\ref{sublemma1}
 that,  for any $a\in\C$,  this factors to give a multiplication map
 of $(\DD^G_{a-q-p}, \, \DD^G_{a})$-bimodules:
\begin{equation} \label{mult}
 \mu_{p,q} :\  (\DD_{a-q})^{\det^{-p}} \otimes_{\DD^G_{a-q}}
(\DD_{a})^{\det^{-q}}   \longrightarrow
    (\DD_{a})^{\det^{-(p+q)}}.
     \end{equation}
 
The aim of this section is to thoroughly understand this: for the appropriate values 
of $c$, we will show that 
$\DD_{c-q}^{\det^{-p}}\cong \QQQ{c-p-q-1}{c-q-1}$ as a $(\UU_{c-p-q},\,\UU_{ c-q-1})$-bimodule
and that $\mu_{p,q}$ is a filtered  isomorphism, thereby proving Theorem~\ref{iso}. Since 
$ \gr  \DD_{c-q}^{\det^{-p}} \cong A^p$, this will also prove Theorem~\ref{main-thm}(1), from 
which the rest of that theorem   follows easily.
As might be expected, the proof is by induction on $p,q$ with the results of the last section providing the
starting point.

\subsection{}  Most modules considered in this paper are subfactors of rings of differential operators 
and, unless we say otherwise, they will then be given the differential operator filtration induced from that 
ambient space.  One exception is with tensor products.
We recall that if $R=\bigcup R_i$ and $S=\bigcup S_j$ are filtered modules over 
some ring $U$ then \emph{the tensor product  filtration}  $\Lambda$  on 
 $R\otimes_US$ is given by $\Lambda_n(R\otimes_US) = \sum_{i+j=n} R_i\otimes S_j$.

   \begin{lemma}\label{products}  Let $p,q\in 
\mathbb{N}$, set $r=p+q$,  and pick $a\in\C$  such that  $a-2, a-3,\dots, a-r$ are all good
{\rm(}this is automatic if $r=0$ or $1${\rm)}.
\begin{enumerate}
\item If $r>0$ then either $a-1$ or $a-r-1$ is good. In the former case $\DD_{a}^{\det^{-r}}$ is a projective 
 right $\DD_{a}^G$-module  while in the latter  it is a projective left $\DD_{a-r}^G$-module.
\item With the tensor product filtration on the left hand side of \eqref{mult}, 
the map $\mu_{p,q}$    is a filtered $(\DD_{a-p-q}^G, \DD_{a}^G)$-bimodule isomorphism. 
\end{enumerate}
 \end{lemma}
 
\begin{proof}   We prove the  two parts of  the  lemma
 by simultaneous induction on $r$. Note that   (1) is vacuous  if $r=0$  and that (2) is 
automatic
   if either  $p=0$ or $q=0$.
Suppose that  $r=p+q = 1$. Then  Proposition~\ref{isom=1} implies that $\DD_{a}^{\det^{-1}}\cong \QQQ{a-2}{a-1}$  
as $(\UU_{a-2},\, \UU_a)$-bimodules    and so    (1)  is given by Theorem~\ref{the trick}(3).

 We may now assume that $0< p,q < r$. This ensures
 that  $a-q-1$ is good and so 
 the  induction hypothesis   implies that    both $\DD_{a-q}^{\det^{-p}} $ and $
 \DD_{a}^{\det^{-q}}$ are projective modules over $\DD_{a-q}^G$.  Moreover, 
  as  one of   $a-r-1$ and 
$a-1$ is good, one of these modules is also a    projective  module on the other side. 
  Thus, once the proof of part~(2) is complete, part~(1) will follow  from 
    Lemma~\ref{proj-tensor}.
   
It remains to prove part~(2).
We first claim that $  \DD_{a-q}^{\det^{-p}} \otimes_{\DD_{a-q}^G}
\DD_{a}^{\det^{-q}}$  is a rank one, torsion-free right $\DD_{a}^G$-module. 
To see this, first   take the differential operator filtration on some 
 $\DD_{b}^{\det^{-s}}$. By  Lemma~\ref{sign-reps}(1),
 $\gr  \DD_{b} ^{\det^{-s}}  \cong A^s$, which  is  torsion-free of  rank one  as a  module over
 $A^0=\C[\h\times\h^*]^W = \gr \DD_{b}^G$. It follows that $\DD_{b}^{\det^{-s}}$ is a rank one torsion-free 
  right $\DD_b^G$-module;  in particular $  \DD_{a-q}^{\det^{-p}} $ and $
\DD_{a}^{\det^{-q}}$  are rank one torsion-free modules on the right.    
By the previous paragraph they are also both projective modules over 
$\DD_{a-q}^G$. Hence, as a right $\DD_a^G$-module,
 $  \DD_{a-q}^{\det^{-p}} \otimes_{\DD_{a-q}^G}\DD_{a}^{\det^{-q}}$ embeds into 
 $  \DD_{a-q}^G\otimes_{\DD_{a-q}^G}\DD_{a}^{\det^{-q}} = \DD_{a}^{\det^{-q}}$, which is a 
 rank one torsion-free right $\DD_{a}^G$-module.  This proves the claim. 

We return to the proof of part~(2).
  It is clear from Lemma~\ref{sublemma1} and Theorem~\ref{first-identity} 
 that $\mu_{p,q}$ is a $(\DD_{a-p-q}^G, \DD_{a}^G)$-bimodule homomorphism 
 and so  we just need to prove that  it  is a filtered  vector space  isomorphism.
 We first show that $\mu_{p,q}$  injective.  
 To see this note that, as above,  $\DD_{a}^{\det^{-r}}$ is a torsion-free 
  right $\DD_a^G$-module.   
Also, $\mathrm{Im}(\mu_{p,q})\neq 0$ since  
 $\mu_{p,q}(\svol^p\otimes \svol^q) = \svol^{p+q} \neq 0$.  Hence $\mathrm{Im}(\mu_{p,q})$
 is a non-zero torsion-free right $\DD_a^G$-module. Since $   \DD_{a-q}^{\det^{-p}} \otimes_{\DD^G_{a-q}}
\DD_{a}^{\det^{-q}}$  is a rank one torsion-free right $\DD_{a}^G$-module, this forces 
 $\mu_{p,q}$ to  be  injective. 

We next prove that $\mu_{p,q}$ is surjective. We  can now  identify 
$  \DD_{a-q}^{\det^{-p}} \otimes_{\DD_{a-q}^G}  \DD_{a}^{\det^{-q}}$ with  its image    
$  \DD_{a-q}^{\det^{-p}} \cdot   \DD_{a}^{\det^{-q}}$ under $\mu_{p,q}$. This is   filtered by the 
image of the 
tensor product filtration and so  \cite[Lemma~6.7(1)]{GS} implies that 
$$\gr \bigl(\DD_{a-q}^{\det^{-p}}\bigr) \cdot\  \  \gr \bigl(\DD_{a}^{\det^{-q}} \bigr)  \
\subseteq \  \gr  \bigl(\DD_{a}^{\det^{-r}}\bigr).$$
By Lemma~\ref{sign-reps}(1) we can regard this multiplication as taking place inside $\C[\h\times\h^*]$ 
in which 
case two further applications of that lemma give
\begin{equation}\label{grsur}
  \bigl( \gr \DD_{a-q}^{\det^{-p}}\bigr) \cdot \
\bigl( \gr \DD_{a}^{\det^{-q}} \bigr) \ 
= \  A^p\cdot A^q = A^{p+q} =  \gr \DD_{a}^{\det^{-r}}.
\end{equation}
 Thus $\mu_{p,q}$ is graded surjective and hence surjective. Since it is immediate that $\mu_{p,q}$
  is a filtered homomorphism, this also proves that it is a filtered isomorphism. This completes the 
  proof of part~(2) and  hence  of the lemma. 
   \end{proof}

 \subsection{}   We can now transfer results from the  $\gr \DD_{c}^{\det^{-p}}$ to the 
 $\QQQ{c-p-1}{c-1}$ and complete the proofs of Theorems~\ref{iso} and~\ref{main-thm} for the modules 
$\QQQ{a}{b}$.

\begin{theorem} \label{warmup} Fix an integer $m\geq 0$ and $c\in \C$ such that the numbers
 $c-1,c-2,\dots,c-m+1$ are good {\rm(}this is automatic if $m=0,1${\rm)}. 
\begin{enumerate}
\item   Under the differential operator
 filtration on both sides there is a filtered isomorphism
 $\Theta_{c,m} : \DD_{c+1}^{\det^{-m}} \iso   \QQQ{c-m}{c}$ of  $(\UU_{c-m},\,\Uc)$-bimodules.
  \item  Under the differential operator
 filtration $\Gamma$   we have $\gr_\Gamma  \QQQ{c-m}{c} = \delta^{-m}A^m e$.
 \end{enumerate}
  \end{theorem}

 \begin{proof} When $m=0,1$, the result follows from Theorem~\ref{first-identity}, respectively
  Corollary~\ref{iso2} and Lemma~\ref{sign-reps}. So we can assume that $m\geq 2$.
 
Consider the following chain of maps:
\begin{eqnarray}  \label{bigchain}
\Theta: \DD_{c+1}^{\det^{-m}} &\stackrel{\alpha}\longrightarrow &  
\DD_{c-m+2}^{\det^{-1}} \otimes_{\DD_{c-m+2}^G} \cdots
\otimes_{\DD_{c}^G}  \DD_{c+1}^{\det^{-1}} \\ 
 &\stackrel{\beta}\longrightarrow
&(\QQQ{c-m}{c-m+1})\otimes_{\UU_{c-m+1}} \cdots \otimes_{\UU_{c-1}} (\QQQ{c-1}{c})\
\stackrel{\gamma}\longrightarrow\ \QQQ{c-m}{c}.  \notag
\end{eqnarray}
Here $\alpha$ is the isomorphism given by iterating the $\mu_{pq}^{-1}$ 
for appropriate $p,q$ and applying Lemma~\ref{products}(2). Similarly, $\beta$ is the map 
$\Psi_{c-m+1}\otimes\cdots \otimes \Psi_{c}$, which is  an isomorphism by 
Proposition~\ref{isom=1} and induction. Finally,  by \eqref{trick-equ} $\gamma$ is an isomorphism induced  by 
multiplication in $\ureg$.   By those same references,  each of these maps is
 a $(\UU_{c-m},\,\Uc)$-bimodule map and hence $\Theta$ is a $(\UU_{c-m}, \, \Uc)$-bimodule isomorphism.

We claim that $\Theta$ is a filtered isomorphism, where 
 the  domain and codomain  are given the filtrations $\Gamma$  induced by the differential operator 
filtration
 in $\de$ and $\D(\hr)\ast W$ respectively.  On the two middle terms we will take the tensor product filtration
 induced from the differential operator filtration on the individual tensor-summands.  

In order to prove this claim, we study the individual maps. 
   By Lemma~\ref{products}(2) each   map $  \mu_{p,q}$  is a  filtered isomorphism
   and hence so is $\mu_{p,q}^{-1}$.  By induction this implies that $\alpha$ is also a filtered 
isomorphism.
   If we take  the tensor product filtration on both 
$$\DD_{c-m+2}^{\det^{-1}} \otimes \cdots \otimes  \DD_{c+1}^{\det^{-1}}
\qquad\text{and}\qquad \QQQ{c-m}{c-m+1}\otimes \cdots \otimes \QQQ{c-1}{c},$$
then  Corollary~\ref{iso2} implies that $\beta$ is then a filtered
isomorphism. 
Next, it is clear  that $\gamma$ is a map of filtered  modules since
 \begin{eqnarray*}    
  \gamma \big(\Gamma_{r_1}(\QQQ{c-m}{c-m+1})\otimes \cdots  \otimes \Gamma_{r_m}
(\QQQ{c-1}{c})\big) &  \\
  \noalign{\vskip 5pt}
  \ = \Gamma_{r_1}(\QQQ{c-m}{c-m+1}) \cdots 
\Gamma_{r_m} (\QQQ{c-1}{c})  \notag \ & \subseteq \ \  \Gamma_{r_1+\cdots + r_m} (\QQQ{c-m}{c}),
\notag
 \end{eqnarray*}
  for all $r_i\geq 0$. However, we do not yet know that $\gamma $ is a filtered isomorphism.
  
  Putting these facts together implies that $\Theta$ is a filtered homomorphism.  
     We now follow    the argument of Corollary~\ref{iso2} to deduce that
   $\Theta$ is a filtered isomorphism. Suppose this is not true. Then, as $\Theta$ is an isomorphism of 
unfiltered objects, 
    there must then exist some $ d\in D=\DD_{c}^{\det^{-m}}$ for which the order of 
    $\Theta(d)$ is strictly less than the order of $d$.
 Hence  $ \text{gr}(\Theta) : \text{gr}_\Gamma D \to \text{gr}_\Gamma ( \QQQ{c-m}{c})$ will have a
nonzero kernel. But we know from Lemma~\ref{sign-reps}(1)
 that  $ \text{gr}_\Gamma D$ is a rank one torsion-free $\C[\h\times \h^*]^W$-module, and so  its image
 must therefore be torsion. Since $\text{gr}_\Gamma (\QQQ{c-m}{c})$ is a submodule of the torsion-free 
$A^0$-module
 $\text{gr}_\Gamma \ureg=\C[\hreg\times \h^*]^We$ this forces $\text{gr}(\Theta)=0$.
 This contradicts the fact that, for example,
 $0\not=e = \Theta(\svol^m)  \in \Theta(\Gamma_0D) \subseteq \Gamma_0(\QQQ{c-m}{c}),$
 where $\svol$ is the function defined in \eqref{sdefn}. This  proves part~(1) of the theorem.
  
 It now follows from  Lemma~\ref{sign-reps}  that 
$\gr_\Gamma(\QQQ{c-m}{c}) \cong \gr_\Gamma \DD_{c+1}^{\det^{-m}} \cong A^{m}\cong \delta^{-m}
A^me$ 
as graded $A^0$-modules, thereby proving that part~(2) holds up to isomorphism.
However, the theorem asserts that there is an equality of graded modules
  $\gr_\Gamma(\QQQ{c-m}{c})  =  e\delta^{-m}A^m$ as subspaces of 
  $\gr_\Gamma \ureg = e\C[\hreg\times\h^*]^W$. This requires a little more work.

  We prove this equality by induction. As in the proof of \cite[Lemma~6.9(2)]{GS}, 
  and for any   $a\in \C$, 
   $$\gr_\Gamma(\QQQ{a-1}{a}) 
=\gr_\Gamma (e\delta^{-1} H_{a}e) = e\delta^{-1}\C[\h\times\h^*]\ast W e
= e\delta^{-1} \C[\h\times\h^*]^{\sgn}e = e\delta^{-1}A^1,$$
  and so the result holds for $m=1$.
Since we have proved that both $\alpha$ and $\beta$ are filtered isomorphisms, 
the fact that $\Theta $ is a filtered isomorphism also implies  that $\gamma$ is a filtered
isomorphism. This can be tautologically reformulated as the statement:
\begin{equation}\label{tautological}
\begin{array}{cc}
&\text{\it the differential operator  and  tensor product filtrations are equal  on}\\
\noalign{\vskip 5pt}
&  \QQQ{c-m}{c}  = (\QQQ{c-m}{c-m+1})\cdot(\QQQ{c-m+1}{c-m+2}) \cdots (\QQQ{c-1}{c}).
\end{array}
\end{equation}
By \cite[Lemma~6.7(2)]{GS} and induction, the multiplication map therefore  induces  a surjection 
$$ \chi: \ e \delta^{-m}A^m \ = \ \gr_\Gamma(\QQQ{c-m}{c-m+1}) \cdots  \gr_\Gamma(\QQQ{c-1}{c})  \ 
\ontoo \ \gr_\Gamma (\QQQ{c-m}{c}).$$  As both sides of this equation
are non-zero  rank one torsion-free $A^0$-modules, $\chi$ must be an equality.
\end{proof}

\subsection{}  We can use the theorem to improve the results from Lemma~\ref{products}.

\begin{corollary}\label{cor6}  Fix an integer $m>0$ 
and  $c\in \C$ such that  $c-1,c-2,\dots,c-m+1$ are all good. Then, $\DD_{c+1}^{\det^{-m}}$
is reflexive as both a left $\DD_{c-m+1}^G$-module and a right $\DD_{c+1}^G$-module
and is the unique such reflexive bisubmodule of $\ureg$.
\end{corollary}

\begin{proof} Combine Theorems~\ref{warmup}  and \ref{the trick}(2).
\end{proof}

\subsection{Remark} \label{P-remark}
Theorem~\ref{warmup} completes the proof of Theorem~\ref{iso} and  Theorem~\ref{main-thm}(1).

\medskip In order to complete the proof of  Theorem~\ref{main-thm} we  need to understand 
the associated graded modules of the $\PPP{c+m}{c}$.
  We expect, but have not pursued, an isomorphism analogous to Theorem~\ref{warmup}(1)
 between $\PPP{c+m}{c}$ and $[\D(\GG)/I_{c+m+1}\D(\GG)]^{\det^m}$.  The proof will need to be a
  little more involved  since the radial component map $\mathfrak{R}$ from \eqref{Rad-defn} 
 actually induces the zero map on  $\D(\GG)^{\det}$. 
This can presumably be circumvented by using
  a ``Fourier transform'' of $\mathfrak{R}$.  We will, however,  take an alternative approach by showing in the
  next   lemma that there   is an easy direct  way to move between the $\mathsf{Q}$'s and $\mathsf{P}$'s that 
  makes such a result unnecessary.

\subsection{}  By   \cite[Remark~2.2]{De} there is an   isomorphism  $\phi_c: \Hc\to {\sf{H}}_{-c}$  
defined by $\phi_c (x) = x$, $\phi_c (D_c(y)) = D_{-c}(y)$ and
$\phi_c(w) = \sgn (w) w$, for $x\in \C[\h]$, $y\in \C[\h^*]$ and $w\in W$.
Note that $\phi_c(e)=e_-$  in the notation from \eqref{cherednik-defn}.
  Since   $\phi_c (\delta) = \delta$  the map $\phi_c$  extends to an
automorphism of $\Hreg= \Hc[\delta^{-1}]$ which will still be written $\phi_c$. The reader should be 
warned  that the action of $\phi_c$ on $\Hreg$ does depend upon~$c$.

\begin{lemma}\label{twist2} Fix $c\in \C$. Then:
  \begin{enumerate}
\item $\phi_{c+1}(z) = \delta^{-2}\phi_c(z)\delta^2$ for all $z\in \ureg$.
\item $\phi_{c+1}(\Uc) = \delta^{-1} {\sf U}_{-c-1} \delta$.
\item
 $\phi_{c+1}(\PPP{c+j}{c}) = \delta^{2j-1}(\QQQ{-c-j-1}{-c-1})\delta$ for  all $j>0$.
 \end{enumerate}
\end{lemma}

\proof  (1)
 By \cite[Proposition~4.9 and (11.33)]{EG},  $\Uc$ is generated as an algebra by 
 $\C[\h]^W$ and $\nabla_{c}^2 = \sum_{i=1}^n D_c(y_i)^2$; thus  $\D(\hr)^W$ is generated 
 by $\C[\hr]^W$ and $\nabla_{c}^2 $. 
     Thus we only need to confirm (1) for elements $pe\in \C[\hr]^We$ 
   and for $\nabla_c^2e$. The former  is obvious since  
      $$\delta^{-2}\phi_c(pe) \delta^2 = \delta^{-2} pe_- \delta^{2} =   pe_- = \phi_{c+1}(pe)
      \qquad\text{for}\ p\in \C[\hreg]^W.$$ 
   By \cite[Theorem~3.1]{heck} we have $\delta^{-1} \nabla_{c+1}^2e_- \delta = \nabla_c^2e$ for all $c\in 
\C$
   (note that our scalar $c$ is the scalar $-k$ in \cite{heck}).
 Thus, (1)
follows from the chain of equalities
   \begin{eqnarray*}\delta^{-2}\phi_c(\nabla_c^2e) \delta^2  = 
   \delta^{-2} \nabla_{-c}^2e_-\delta^2 \ = \  \delta^{-1} \nabla_{-c-1}^2e \delta =
    \phi_{c+1}(\delta^{-1} \nabla_{c+1}^2e_- \delta) = \phi_{c+1} ( \nabla_c^2e).\end{eqnarray*}

(2)  This is equivalent to the assertion  that $\phi_{-c-1}({\sf U}
_{-c-1}) = \delta \Uc\delta^{-1}$. By \cite[Proposition~4.9 and (11.33)]{EG}, again, it is enough to confirm 
this for $\C[\h]^We$ and $\nabla_{-c-1}^2e$. This is trivial for $\C[\h]^We$ while \cite[Theorem~3.1]{heck}
gives $\phi_{-c-1}(\nabla_{-c-1}^2e) = \nabla_{c+1}^2e_- = \delta \nabla_c^2e \delta^{-1},$ 
as required. 

(3) For  $j=1$ we
have $$\phi_{c+1}(\PPP{c+1}{c}) = \phi_{c+1}(e{\mathsf H}_{c+1}\delta e)
= e_-{\mathsf H}_{-c-1}\delta e_- = \delta e \delta^{-1} {\mathsf
H}_{-c-1} e \delta = \delta (\QQQ{-c-2}{-c-1})
\delta,$$
as desired. Now assume that the result holds for some $j\geq 1$. By part~(1),
$\phi_{c+j+1}(z) = \delta^{-2j} \phi_{c+1} (z) \delta^{2j}$
for $z\in \ureg$. Applying  the case $j=1$ gives $$\delta (
\QQQ{-c-j-2}{-c-j-1})\delta = \phi_{c+j+1}(\PPP{c+j+1}{c+j}) =
\delta^{-2j}(\phi_{c+1}(\PPP{c+j+1}{c+j}))\delta^{2j}.$$ Equivalently
$$\phi_{c+1}(\PPP{c+j+1}{c+j}) =
\delta^{2j+1}(\QQQ{-c-j-2}{-c-j-1})\delta^{1-2j}.$$ To complete the proof, by induction, we
obtain \begin{align*} \phi_{c+1}(\PPP{c+j+1}{c}) &=
\phi_{c+1}(\PPP{c+j+1}{c+j})\phi_{c+1}(\PPP{c+j}{c}) \\
  &=
\delta^{2j+1}(\QQQ{-c-j-2}{-c-j-1})\delta^{1-2j}\cdot\delta^{2j-1}
(\QQQ{-c-j-1}{-c-1})\delta \\ 
    &= \delta^{2j+1} (\QQQ{-c-j-2}{-c-1})\delta.
\tag*{$\Box$}\end{align*}

\subsection{Completion of the proof
of Theorem~\ref{main-thm}.}  By Remark~\ref{P-remark} it only remains to prove Theorem~\ref{main-thm}(2).
Changing notation slightly, we consider $\PPP{d+m}{d}$ where  $d\in\C$ is chosen so 
 that the numbers $d+1,d+2, \dots, d+m-1$ are all good. Thus we need to prove that 
 $\gr \PPP{d+m}{d} = A^m\delta^me$.  
 
 Consider $\phi_{d+1}$. By construction, this morphism   preserves the differential operator
 filtration on $\Hreg$ and   $\gr \phi_{d+1}$ is   then the automorphism that is the identity on
  $\C[\h\times \h^*]$ and sends $w\in W$ to $\sgn(w)w$. Moreover,  if  $c=-d-1$ then 
$c-1,c-2,\dots, c-m+1$ are also good. 
 Thus Lemma~\ref{twist2} and Theorem~\ref{warmup} can be applied to show that 
  \begin{eqnarray*}  (\gr\phi_{d+1})( \gr \PPP{d+m}{d}) & = &  \gr\bigl( \phi_{d+1}(\PPP{d+m}{d}) \bigr) \ = \ 
  \gr\bigl(\delta^{2m-1} (\QQQ{-d-m-1}{-d-1}) \delta\bigr) \\ 
  \noalign{\vskip 7pt}
  & = & \gr\bigl(\delta^{2m-1} (\QQQ{c-m}{c}) \delta\bigr) \ = \    \delta^{2m-1} (\delta^{-m}A^me)\delta 
\ =\ \delta^{m-1} A^m  \delta  e_- 
\end{eqnarray*}
and hence $  \gr(\PPP{d+m}{d}) =  \gr \phi_{d+1}^{-1}\bigl( \delta^{m-1} A^m  \delta  e_-\bigr)
=  \delta^m A^m e,$ as required. \qed

\subsection{} \label{shiftrelation} We finish this section by making explicit the connection between 
shifting by ${\sf Q}$'s, by ${\sf P}$'s and applying the various morphisms $\phi_c$'s of 
Section~\ref{twist2}. The equality $\phi_{c+1}(\Uc) = \delta^{-1}{\sf U}_{-c-1}\delta$ from 
Lemma~\ref{twist2}(2)  induces an equivalence of categories 
$$\Omega_c : \Lmod{\Uc}\stackrel {\sim} \longrightarrow \Lmod{{\sf U}_{-c-1}}, \qquad V\mapsto V,$$ 
where $V\in   \Lmod{\Uc}$ becomes a ${\sf U}_{-c-1}$-module via $z\ast v = \phi_{-c-1}(\delta^{-1} z \delta)v$ 
for  $z\in {\sf U}_{-c-1}$ and $v\in V$.

\begin{proposition} Let $c\in \C$ and set $d = - c -1$. For any integer $m\geq 0$ there is an isomorphism 
of functors $\PPP{d+m}{d}\otimes_{{\sf U}_d}\Omega_c(-) \cong \Omega_{c-m}\circ (\QQQ{c-m}{c}
\otimes_{\Uc} (-))$ that makes the following diagram commute 
$$ \begin{CD} {\Lmod{\Uc}} @> { \ \ \QQQ{c-m}{c}\otimes_{\Uc} (-)  \ \  }>> \Lmod{{\sf U}_{c-m}} \\ @V{ \Omega_c}
VV @VV\Omega_{c-m}V \\
\Lmod{{\sf U}_{d}} @>{\ \ \PPP{d+m}{d}\otimes_{{\sf U}_d  } (-) \ \  }>> \Lmod{{\sf U}_{d +m}} \end{CD}$$
\end{proposition}

\begin{proof}
Let $V$ be a $\Uc$-module. By Lemma~\ref{twist2}(3) $\QQQ{c-m}{c} = \delta^{1-2m}\phi_{-c}(\PPP{d+m}
{d})\delta^{-1}$, so we are asserting the existence a natural isomorphism between $\PPP{d+m}{d}
\otimes_{{\sf U}_d} V$ (with ${\sf U}_d$-action on $V$ induced from $\Omega_c$) and $\delta^{1-2m}
\phi_{-c}(\PPP{d+m}{d})\delta^{-1}\otimes_{\Uc} V$ (with ${\sf U}_{d+m}$-action induced from $\Omega_
{c-m}$). The only choice is the mapping $p\otimes v \mapsto \delta^{1-2m}\phi_{-c}(p)\delta^{-1}\otimes_
{\Uc} v$ for $p\in \PPP{d+m}{d}$ and $v\in V$. We need to check that this is well-defined and that it is a $
{\sf U}_{d+m}$-module homomorphism.

Pick $z\in {\sf U}_d$, $p\in \PPP{d+m}{d}$ and $v\in V$. For well-definedness we have \begin{eqnarray*}
pz\otimes v - p\otimes z\ast v &\!\mapsto\!&   
\delta^{1-2m}(\phi_{-c}(pz)\delta^{-1} \otimes_{\Uc} v - \phi_{-c}
(p)\delta^{-1}\otimes_{\Uc} z\ast v) \\ & \!=\! &
 \delta^{1-2m}\phi_{-c}(p)\left(\phi_{-c}(z)\delta^{-1} \otimes_
{\Uc} v -\delta^{-1}\otimes_{\Uc} \phi_{-c-1}(\delta^{-1}z\delta)
 v \right)\\ & \!=\!&  \delta^{1-2m}\phi_{-c}(p)\left
(\delta^{-2}\phi_{-c-1}(z)\delta \otimes_{\Uc} v -\delta^{-1}\otimes_{\Uc} \phi_{-c-1}(\delta^{-1}z\delta) v 
\right) = 0, 
\end{eqnarray*}
where in the first equality of the last line we used Lemma~\ref{twist2}(1).

Now let $z\in {\sf U}_{d+m}$. Using Lemma~\ref{twist2}(1) and induction we have
\begin{eqnarray*}
zp\otimes v &\mapsto & \delta^{1-2m}\phi_{-c}(zp)\delta^{-1}\otimes_{\Uc} v  \\ & = & \delta^{-2(m-1)}\phi_
{-c}(\delta^{-1}z\delta ) \delta^{2(m-1)}\delta^{1-2m}\phi_{-c}(p)\delta^{-1}\otimes_{\Uc} v \\ & = & \phi_{-c
+m-1}(\delta^{-1}z\delta) \delta^{1-2m}\phi_{-c}(p)\delta^{-1}\otimes_{\Uc} v \\ &= & z\ast \left(\delta^
{1-2m}\phi_{-c}(p)\delta^{-1}\otimes_{\Uc} v \right).
\end{eqnarray*}
This confirms ${\sf U}_{d+m}$-equivariance.
\end{proof}


\section{Shift functors for $\D$-modules and Cherednik algebras}\label{shift-section}

\subsection{}  The morphism $\PPP{c+m}{c}\otimes - : \UU_c\operatorname{-mod}
\to \UU_{c+m}\operatorname{-mod}$  is fundamental to the representation theory of $\UU_c$, as is illustrated
by much of \cite{GS2}.  There is a similar translation functor for twisted $\D$-modules on projective space  
given by tensoring with $\mathcal{O}(mn)$ (see Section~\ref{commute-sect} for the precise definition).
 In this section we show that these functors are naturally intertwined by  hamiltonian reduction, thereby 
 proving Theorem~\ref{shiftthm} from the introduction.   Before stating the result we will need some definitions.

\subsection{Projectivisation.}\label{mor1} The quantum hamiltonian
reduction of Theorem~\ref{first-identity} can be performed in two steps: first with respect
to the  subgroup   of scalar matrices $\C^\tstar\subset G=GL(V)$; and
then with respect to the subgroup $SL(V)\subset G$. 
Recall that $G$ acts diagonally on  $\GG=\g\times V$ and so 
 $\C^\tstar$ acts trivially on $\g$ and by dilations $\lambda\circ v=\lambda v$  on $V=\C^n$.
 As in \eqref{hamilton-defn},   the identity matrix in $\gl(V)$ is written $\eu$.   Put 
$$\Vo =  V\sminus\{0\}; \quad \P=\P(V)=\Vo\!/\C^\tstar;\quad \GG^\circ=
\g\times\Vo;\quad\text{and}\quad\X=\g\times\P.
$$

\begin{lemma} \label{mor1-lemma}
Let $c\in\C$ and assume that if $n=2$ then $2c\notin \mathbb{Z}_{\leq 0}$. Then restriction provides a 
natural isomorphism
\begin{equation}\label{mor1-equ}
\chi : \D(\GG)/\D(\GG)\tau(\eu- nc)\stackrel{\sim}\longrightarrow 
\Gamma\big(\GG^\circ,\,\D_{\GG^\circ}/\D_{\GG^\circ}\tau(\eu-nc)\big).
\end{equation}
\end{lemma}

\begin{proof} Without loss of generality we need only consider $V$ instead of $\GG$.

Set $v = \tau(\eu -nc)$. Multiplication by $v$ on the right
yields a short exact  sequence of sheaves on $V$:
$$0\to \D_V \to\D_V \to\D_V/\D_Vv\to 0.$$
Hence on $V^\circ$ we obtain the long exact sequence
$$\begin{array}{rl}
0 \too\Gamma (V^\circ, \D_V)\ \stackrel{\times v}\too & \Gamma (V^\circ, \D_V)
 \too\Gamma (V^\circ, \D_V/\D_Vv)\too \\
\noalign{\vskip 7pt}
& \too H^1(V^\circ,\D_V)\ \stackrel{\times v}\too\ H^1(V^\circ , \D_V) \too \ldots
\end{array}$$
 Assume first that $n\geq 3$. Then 
$H^1(V^\circ, \mathcal{O}_V) \cong H^2_{\{0\}}(V, \mathcal{O}_V) = 0$ since 
$H^i_{\{0\}}(V, \mathcal{O}_V) =0 $ for all $i<n=\dim V$.  
Since $\D_V$ is free  as a sheaf of 
${\mathcal{O}}_V$-modules, it follows that $H^1(V^\circ,\D_V)=0$. Since $\C[V]$ is just a polynomial ring
and $\dim V\smallsetminus V^\circ \leq \dim V-2$, we have 
$\C[V^\circ] = \Gamma (V^\circ, \mathcal{O}_{V^\circ}) = \C[V]$ by Hartog's theorem. By freeness,  
$\Gamma (V^\circ, \D_{V^\circ}) = \D(V)$ and so we are done if $n\geq 3$.

Now assume that $n=2$ and write  $\C[V] = \C[x,y]$. A simple calculation in \v{C}ech cohomology using 
the open sets $D(x) = V\setminus \{ x=0 \}$ and $D(y)= V\setminus \{ y=0\}$ shows that 
$$H^1(V^\circ , \mathcal{O}_V) \cong \frac{\C[x^{\pm 1}, y^{\pm 1}]}{(\C[x^{\pm 1}, y] + \C[x, y^{\pm1}])}.$$ 
Denote this space by $S$, so that $H^1(V^\circ, \D_V) = \bigoplus_{i,j\geq 0} S\partial_x^i\partial_y^j,$ 
where $\partial_x=\partial/\partial x$, etc.

We claim that the mapping $H^1(V^\circ,\D_V)\stackrel{\times v}\to H^1(V^\circ , \D_V)$ is injective if 
$2c\notin \Z_{\geq 0}$. To see this let 
$$0\not= z = \sum_{i,j\geq 0} s_{i,j} \partial_x^i\partial_y^j
\in H^1(V^\circ,\D_V),$$
 for some $s_{i,j}\in S.$
By \eqref{hamilton-defn2} we have  $\tau(\eu) = -x\partial_x - y\partial_y$ and~so
\begin{eqnarray*} 
zv \ = \ z\tau(\eu - 2c) &=& \sum_{i,j\geq 0} \left( s_{i,j} \tau(\eu) - (i+j +2c)s_{i,j}\right)\partial_x^i\partial_y^j 
\\  \noalign{\vskip 5pt}
& = & \sum_{i,j \geq 0} \left(s_{i-1,j}x + s_{i,j-1}y - (i+j + 2c)s_{i,j}\right) \partial_x^i\partial_y^j.\end{eqnarray*}
 Thus if $zv = 0$ then $s_{i-1,j}x + s_{i,j-1}y - (i+j+2c)s_{i,j} =0$ for all $i,j\geq 0$. But consideration of the
  lowest degree $(i_0,j_0)$ of nonzero monomials in $z$   then shows that $i_0+j_0+2c = 0$, 
  contradicting the hypothesis that 
  $2c \notin \Z_{\leq 0}$. It follows that the morphism $$\Gamma(V, \D_V) \cong \Gamma (V^\circ, \D_V)
\to\Gamma (V^\circ, \D_V/\D_Vv)$$ is surjective, as required.
\end{proof}

\noindent
{\bf Remark.} The lemma fails when $n=2$ and $2c \in \mathbb{Z}_{\leq 0}$. In the notation of the above  
proof we have the following equality in $\D_V(D(xy))$:
 \begin{eqnarray*} 
 x^{-1}y^{-1}\partial_x^{-2c} v & = &- x^{-1} \partial_x^{-2c}\partial_y - y^{-1}\partial_x^{-2c+1}.
 \end{eqnarray*} 
 Thus $x^{-1} \partial_x^{-2c}\partial_y = - y^{-1}\partial_x^{-2c+1}$ in $(\D_V/\D_Vv)(D(xy))$. It 
 follows that the sections $x^{-1}\partial_x^{-2c}\partial_y \in (\D_V/\D_Vv)(D(x))$ and
  $-y^{-1}\partial_x^{-2c+1}\in (\D_V/\D_Vv)(D(y))$ extend to a section of $(\D_V/\D_Vv)(V^\circ)$ 
  which is not in the image of $\chi$.

 \subsection{}\label{c-red-sect}
Following \cite[Section~5.1]{GG}, for each $c\in\C$ we introduce an algebra  
\begin{equation}\label{c-red}\D_{c}(\X) \ \defn \
\left(\frac{\displaystyle  \D(\GG)}{\displaystyle  \D(\GG)\tau(\eu-nc)}\right)^{\C^\tstar}.
\end{equation}
The algebra on the right hand side  of \eqref{c-red} is a quantum hamiltonian reduction
 with respect to the group $\C^\tstar$ at the point $c$.

 Using Theorem~\ref{first-identity},   we can apply the  formalism of hamiltonian reduction,
as outlined in \cite[Section~7]{GG}, to $\D_{c+1}(\X)$ and $\Uc$.
Let $\Lmod{\bigl(\D_{c+1}(\X), SL(V)\bigr)}$ denote the category whose objects are the finitely generated 
$SL(V)$-equivariant $\D_{c+1}(\X)$-modules on which the action of $\sl(V)$ induced from the 
$SL(V)$-equivariance agrees with the action of $\sl(V)$ induced from the homomorphism 
$\sl(V) \stackrel{\tau}\longrightarrow \D(\GG)^{\C^\tstar} \longrightarrow \D_{c+1}(\X)$. 
We have the following  {\em functor of hamiltonian reduction}:
\beq{ham_fun}
\BH_c:\ \Lmod{(\D_{c+1}(\X), SL(V))} \too \Lmod{\UU_{c}}; \ 
\quad \CF\mapsto \BH_c(\CF)=\CF^{SL(V)}.\eeq

Since $SL(V)$ is reductive, the functor $\BH_c$  is {\em exact} and, by \cite[Proposition 7.1]{GG}, 
it has   a left adjoint 
$$\begin{aligned}
{}^\top\BH_c:\
\Lmod{\UU_{c}} & \too\Lmod{(\D_{c+1}(\X), SL(V))};\quad  \\ \noalign{\vskip 5pt}
M   & \ \mto \  (\D_{c+1}(\X)/\D_{c+1}(\X)\tau(I_{c+1})) \otimes_{\UU_c}\,M.
\end{aligned}
$$

\subsection{} Next, let $\D_{\X,c}$ denote the  sheaf of $(-nc)$-twisted differential
 operators on  $\X$. It follows from \cite[Equation~5.1]{GG} that  $\D_{\X ,c}$    has global sections 
$\Gamma(\X,\D_{\X,c}) = \D_c(\X)$. For any $m\in\Z$, let  $\oo(m)$ be the pull-back of
the standard line bundle $\oo_\P(m)$ via the   projection $\X=\g\times\P\to\P$.
Tensoring with  $\oo(nm)$ yields a functor
\beq{transl}
\Lmod{\D_{\X,c}}\too \Lmod{\D_{\X,c-m}},\quad
{\mathcal F} \mapsto  \oo(nm)\o_{\oo_\X} {\mathcal F}.
\eeq

\subsection{} \label{commute-sect}
 Assume now that $n(c+1)\in\C\smallsetminus \mathbb{Z}_{>0}$.
As in \cite[Proposition~5.4]{GG}, we can apply 
 the Beilinson-Bernstein theorem to give  an equivalence of categories
\begin{equation}\label{BB-equivalence}
\Gamma(\X,-):\
 \Lmod{\D_{\X,c+1}}\iso \Lmod{\D_{c+1}(\X)},\quad
{\mathcal F} \mapsto \Gamma(\X,{\mathcal F}).
\end{equation}
We write
$$
\BS^m:\
\Lmod{\D_{c+1}(\X)}\too\Lmod{\D_{c-m+1}(\X)},
\quad F\mapsto F(nm)
$$
for the functor that corresponds to the functor
\eqref{transl} via the  Beilinson-Bernstein equivalence.

\begin{theorem}\label{commute} 
Fix $c\in \C$ and a positive integer $m$ such that each of the numbers
$c-1, c-2,$ $\ldots , c-m+1$ is good {\rm and } $n(c+1)\notin \mathbb{Z}_{> 0}$
 {\rm (}respectively, $n(c+1)\notin \Z$ if $n=2${\rm )}. Then
there is an isomorphism of functors $\BH_{c-m}\circ \BS^m\circ {}^{\top}\BH_{c}( - ) \cong
\QQQ{c-m}{c}\otimes_{\UU_{c}} ( - )$ that makes
the following diagram commute
$$
\xymatrix{ \strut
{\Lmod{(\D_{c+1}(\X), SL(V))}} \ 
\ar[rr]^<>(0.5){\BS^m}&& \ 
{\Lmod{(\D_{c-m+1}(\X), SL(V))}}\ar[d]^<>(0.5){\BH_{c-m}}\\
\Lmod{\UU_{c}} \ar[u]^<>(0.5){{}^{\top}\BH_{c}} \  
\ar[rr]^<>(0.5){\QQQ{c-m}{c}\otimes_{\UU_{c}}(-)}&& \ \Lmod{\UU_{c-m}}.
}
$$
\end{theorem}
 
\subsection{}\label{commute1}
Before proving the  theorem we note that, by
 Proposition~\ref{shiftrelation}, it has  the following 
equivalent formulation in terms of the ${\sf P}$'s.  For $d\in \C$ we set 
 $$\widetilde{\BH}_d  \ \defn\  \Omega_{-d-1}\circ \BH_{-d-1}: \Lmod
{(\D_{-d}(\X), SL(V))} \longrightarrow \Lmod{\UU_d},$$
where $\Omega_d$ is defined in \eqref{shiftrelation}.
This  has  left adjoint 
$${}^\top\widetilde{\BH}_d = {}^\top\BH_{-d-1}\circ \Omega_{-d-1}^{-1} =   {}^\top\BH_{-d-1}\circ \Omega_d.$$
Then, using the fact that $d$ is good if and only if $-d-1$ is good, we obtain:

\begin{corollary}  Fix
$c\in \C$ and a positive integer $m$ such that   $c+1,\dots, c+m-1$ are all 
good and that $nc\notin \Z_{<0}$ {\rm (}respectively $nc\notin \Z$ if $n=2${\rm )}. 
Then  there is a commutative diagram:
$$
\xymatrix{ \strut
{\Lmod{(\D_{-c}(\X), SL(V))}} \ 
\ar[rr]^<>(0.5){\BS^m}&& \ 
{\Lmod{(\D_{-(c+m)}(\X), SL(V))}}\ar[d]^<>(0.5){\widetilde{\BH}_{c+m}}\\   
\Lmod{\UU_{c}} \ar[u]^<>(0.5){{}^{\top}\widetilde{\BH}_{c}}   
 \  \ar[rr]^<>(0.5){\PPP{c+m}{c}\otimes_{\UU_{c}}(-)}&& \ \Lmod{\UU_{c+m}}
}
$$ 
\vspace{-1.0cm}

 \qed
  \end{corollary} 
  
\noindent{\bf Remark.}  
Theorem~\ref{shiftthm} is a special case of  this corollary.

\subsection{} \label{D-routine}
Recall that $\C^\tstar \subset GL(V)$ is the central subgroup consisting of multiples of the 
identity matrix $\Id$. Given a $GL(V)$-representation $E$ we will denote by $E^{(m)}$ the set of 
semi-invariants $\{ e  \in E: (z\Id)\cdot e = z^{mn}e \text{ for all }z\in \C^\tstar\}$. For $d\in \C$, we define
$$\DDD{d-m}{d} \ \defn \  \left[\frac{\D(\GG)}{\D(\GG)\tau(\eu-nd)}\right]^{(-m)}
\qquad \text{for }m\in\Z.
$$  Note  that  $N^{\det^{-m}} = \left(N^{SL}\right)^{(-m)}$ for any $GL(V)$-module $N$.
By Lemma~\ref{sublemma1},  
 $\DDD{d-m}{d}$ has a natural $(\D_{d-m}(\X),\D_{d}(\X))$-bimodule  structure and it also 
has the following useful properties.

\begin{lemma}  {\rm (1)}  
Fix $c\in\C$  and a positive integer $m$  such that $n(c+1)\not\in \mathbb{Z}_{> 0}$
(respectively $n(c+1)\notin \Z$ if $n=2$). Then, for any $F\in {\Lmod{(\D_{c+1}(\X), SL(V))}},$ 
\begin{equation}\label{D-routine1}
\BS^m(F) \  = \ \Gamma\big(\X, \oo(nm)\o_{\oo_\X}\D_{\X,{c+1}}\big)\o_{\D_{c+1}(\X)}F \ \cong \ 
 \DDD{c-m+1}{c+1} \o_{\D_{c+1}(\X)} { F}
\end{equation}
{\rm (2)}   Assume  that  $c-j$ is good for $1\leq j\leq m-1$.
There are isomorphisms of $(\UU_{c-m},\, \UU_{c})$-bimodules
$$ \left(\frac{\strut\displaystyle\DDD{c-m+1}{c+1}}{\displaystyle\strut
\DDD{c-m+1}{c+1}\tau( I_{c+1})}\right)^{SL}    
\!\!  \cong \  \left(\frac{\displaystyle\strut\D(\GG)}{\displaystyle\strut\D(\GG)\tau(I_{c+1}) }\right)^{\det^{-m}}
\!\! \cong \    \QQQ{c-m}{c}.$$\end{lemma}

 \begin{proof}(1)  The first equality in \eqref{D-routine1}  follows from 
the Beilinson-Bernstein equivalence \eqref{BB-equivalence} and the definition of $\BS^m$. For the displayed
 isomorphism, we consider the principal $\C^{\times}$-bundle $p: V^{\circ} \longrightarrow \mathbb{P}$.
 Then by equivariant descent  
 (see  \cite[Chapter~VII, Section~1]{SGA1}) we have
  \begin{eqnarray*}\mathcal{O}_{\mathbb{P}}(mn) \otimes_{\mathcal{O}_{\mathbb{P}}} \D_{\mathbb{P},c+1} 
  &\cong & p_{\ast}\left( (\mathcal{O}_{V^{\circ}}\otimes_{\C} \det{}^{m})
  \otimes_{\mathcal{O}_{V^{\circ}}} \frac{\D_{V^{\circ}}}{\D_{V^{\circ}}\tau({\eu } - n(c+1))} \right)^{\C^{\times}} \\ 
     \noalign{\vskip 7pt}
  & \cong & p_{\ast} \left(\frac{{\D}_{V^{\circ}}}{\D_{V^{\circ}}\tau({\eu } - n(c+1))}  \right)^{(-m)}.
  \end{eqnarray*} It follows that 
   \begin{eqnarray}\label{D-routine2}
   \Gamma (\mathbb{P}, \mathcal{O}_{\mathbb{P}}(nm) \otimes_{\mathbb{P}} \D_{\mathbb{P},c+1})
  &  \cong & 
   \Gamma \left(\mathbb{P}, \, p_{\ast} \Big(\frac{\D_{V^{\circ}}}{\D_{V^{\circ}}\tau({\eu} - n(c+1))}\Big)\right)^{(-m)}  \\
   \noalign{\vskip 9pt}
 &   \cong &   \Gamma \left(V^{\circ}, \, \, \frac{\D_{V^{\circ} }}{\D_{V^{\circ}}\tau({\eu } - n(c+1)) } \right)^{(-m)}\notag
    \end{eqnarray}
   Now $\D_{\GG^\circ}=\D_\g\otimes \D_{V^\circ}$, with $\C^\times$ acting trivially on $\D_\g$. 
   Therefore,  combining 
  \eqref{mor1-equ} with \eqref{D-routine2} gives the following  isomorphism of  $\big(\D_{c-m+1}(\X),\, \D_{c+1}(\X)\big)$-bimodules:
    \begin{eqnarray*}
    \DDD{c-m+1}{c+1} &\cong & \left(\frac{\D(\GG)}{\D(\GG)\tau(\eu-n(c+1))}\right)^{(-m)}  
\ \cong \     \Gamma\left(\GG^\circ,\, \frac{\D_{\GG^\circ}}{\D_{\GG^\circ}\tau(\eu-n(c+1))}\right)^{(-m)} \\
  \noalign{\vskip 7pt} & \cong & 
    \Gamma (\mathbb{P}, \mathcal{O}_{\mathbb{P}}(nm) \otimes_{\mathbb{P}} \D_{\mathbb{P},c+1})\otimes\D(\g)
   \ \cong \     \Gamma\big(\X, \oo(nm)\o_{\oo_\X}\D_{\X,{c+1}}\big),
    \end{eqnarray*}
as required.

(2) As a morphism of $\D(\GG)^{SL}$-modules, the 
 first isomorphism follows from the observation  that  $N^{\det^{-m}} = \left(N^{SL}\right)^{(-m)}$ for any $GL(V)$-module $N$.
  By Lemma~\ref{sublemma1} this restricts to give an  isomorphism  of $(\D_{c-m+1}(\X),\D_{c+1}(\X))$-bimodules;
equivalently of  $(\UU_{c-m},\, \UU_{c})$-bimodules. The second isomorphism is just   Theorem~\ref{iso}.
 \end{proof}

\subsection{Proof of Theorem~\ref{commute}} Let $M\in \Lmod{\UU_{c}}.$ By Lemma~\ref{D-routine}(1)
\begin{eqnarray*}
 \BS^m\circ{}^\top\BH_c(M) 
 & = & \DDD{c-m+1}{c+1}\otimes_{\D_{c+1}(\X)} 
\left(\frac{\displaystyle\strut  \D_{c+1}(\X)}{ \displaystyle \strut \D_{c+1}(\X)\tau( I_{c+1})} \otimes_{\Uc} M\right)
\\ \noalign{\vskip 8pt}
& \cong& \   \left(\frac{\displaystyle\strut  \DDD{c-m+1}{c+1}}{ \displaystyle \strut \DDD{c+m-1}{c+1
}\tau(I_{c+1})}\right) \otimes_{\Uc} M.
\end{eqnarray*}
Thus Lemma~\ref{D-routine}(2)   implies that 
\begin{equation*}\label{mor_prop}
\begin{array}{rl}
  &  \BH_{c-m}\circ\, \BS^m\circ {}^\top\BH_{c}(M)   \ \cong  \ \left(
  \bigg(\,\,\frac{\displaystyle\strut  \DDD{c-m+1}{c+1}}{ \displaystyle \strut \DDD{c+m-1}{c+1}\tau(I_{c+1})}\,\bigg)
   \otimes_{\Uc} M\right)^{SL}  \\  \noalign{\vskip 10pt}
  &\qquad \qquad  
   = \   \left( \frac{\displaystyle\strut  \DDD{c-m+1}{c+1}}{ \displaystyle
    \strut \DDD{c+m-1}{c+1}\tau(I_{c+1})}\right)^{SL}
\hskip-3pt   \otimes_{\Uc} M       \ \cong \  \QQQ{c-m}{c} \otimes M.
\end{array}
\end{equation*}
\vspace{-1.3cm} 

\qed
\vspace{0.8cm}

\section{Characteristic cycles}\label{char-section}

\subsection{} Let $M$ be a filtered, finitely generated $\Uc$-module.
In \cite{GS2} the authors used the $\mathbb{Z}$-algebra associated to the modules $\{\PPP{c+a}{c}\}$ to
 construct a characteristic cycle $ \ch(M)$ inside the Hilbert scheme $ \Hilb$ that then proved useful in 
 studying the representation theory of $\Uc$.  Using quantum hamiltonian 
 reduction, the authors of \cite{GG} define a second such characteristic cycle.
 This leads to the natural question of whether these varieties are equal;  see \cite[(7.17)]{GG}. In this section we 
 show that this is indeed the case.

\subsection{Hilbert schemes}  \label{hilb-sect}
We write $\Coh(X)$ for the category of coherent sheaves on a scheme $X$. If $B = \bigoplus_
{m\geq 0} B_m$ is a finitely generated graded commutative algebra let $\Lmof{B}$ denote
 the  category  finitely generated  graded left  $B$-modules and   write $\mathbb{F}(M)$ for  the 
coherent sheaf on  the scheme $\Proj B$   corresponding to the module $M\in \Lmof{B}$.

  Set $A=\bigoplus_{m\geq 0}\ A^m,$ where the $A^m$ are defined as
in \eqref{moment}; thus \cite[Proposition~2.6]{Ha2} implies that $\Proj(A) \cong \Hilb$, the  Hilbert scheme of
$n$ points in $\C^2$. Following  \cite{Ha2} and \cite{Na},  there is  the following diagram of schemes 
over $(\h\times\h^*)/ W$:
\beq{hilb}
\xymatrix{
\M\en&\en\M\cyc\en\ar@{_{(}->}[l]_<>(0.5){j}\ar@{->>}[r]^<>(0.5){p}&\en\Hilb
 \ \cong \ 
\Proj A\ \cong\  \Proj\left(\bigoplus_{m\geq 0}\ \C[\M]^{\det^{-m}}\right)
}
\eeq
where  $\M=\mu_\GG^{-1}(0)$
 is  defined as in \eqref{M} with open subvariety   \[ \M\cyc = \bigl\{(X,Y,v,w)\in\M :  \C[X,Y]v=V\bigr\}.\]
In more detail,   $\M\cyc$ is a smooth $GL(V)$-variety and   the map  $j: \M\cyc\into\M$
is a  $GL(V)$-equivariant Zariski open imbedding.
The map $p$ in \eqref{hilb}  is a {\em universal
geometric quotient} morphism that makes
$\M\cyc$ a principal $GL(V)$-bundle over $\Hilb$ (see \cite[Proof of Theorem~1.9]{Na}).
Finally, the penultimate isomorphism in \eqref{hilb}
was proved in \cite[Proposition~2.6]{Ha2} while  the last isomorphism follows from  \eqref{A}.

\subsection{} \label{7.3}
We can also construct $\Hilb$ via 
$$T^*\X = \{ (X,Y,v,w)\in \g\times \g \times V^{\circ} \times V^*: w(v) = 0 \}/\C^{\times},$$ where the action of 
$\C^{\times}$ arises from the scalar matrices in $GL(V)$; thus it acts  only on $V^{\circ}\times V^*$. 
In more detail, since 
$\Tr( [X,Y] + vw) = w(v)$, then, as in \eqref{moment}, 
$$\mu_{\X}^{-1}(0) \defn \{ (X,Y,v,w) \in \g\times \g \times V^{\circ} \times V^*: [X,Y] + vw = 0 \} /\C^{\times}$$
 is a closed subvariety of $T^*\X$. Inside $\mu_{\X}^{-1}(0)$ we have the   the Zariski open subset 
  $$\mu_{\X}^{-1}(0)\cyc \ \defn \  \left\{ (X,Y,v,w)\in  \mu_{\X}^{-1}(0)  : v\in V^\circ \ \text{and}
   \C[X,Y]v = V\right\}\ \subset  \mu_{\X}^{-1}(0) .$$ As before, we have a diagram schemes
\begin{equation*}
\xymatrix{
\mu_{\X}^{-1}(0)\en&\en\mu_{\X}^{-1}(0)\cyc\en\ar@{_{(}->}[l]_<>(0.5){j_{\X}}\ar@{->>}[r]^<>(0.5){p_{\X}}&\en\Hilb
}
\end{equation*}
where $p_{\X}$ is a universal geometric quotient making a principal $PGL(V)$-bundle.
There is, moreover, a commutative diagram 
\begin{equation}\label{hilbcomm}
\xymatrix{
\M\en&\en\M^\circ\en\ar@{_{(}->}[l]_<>(0.5){j_1} \ar@{->>}[d]_{\iota} &
 \en\M\cyc\en\ar@{_{(}->}[l]_<>(0.5){j_0}\ar@{->>}[dr]^<>(0.5){p} \ar@{->>}[d]_{\iota\cyc} \ar@(ul,ur)[ll]_{j}& \\ 
& \mu_{\X}^{-1}(0)\en&\en\mu_{\X}^{-1}(0)\cyc\en\ar@{_{(}->}[l]_<>(0.5){j_{\X}}\ar@{->>}[r]^<>(0.5){p_{\X}}&
 \en\Hilb 
}
\end{equation}  
 where $\mathcal{M}^\circ = \{ (X,Y,v,w)\in \mathcal{M} : v\neq 0\}$
 and the vertical maps $\iota, \iota\cyc$ are principal $\C^{\times}$-bundles.
 
\subsection{} 
We will require a special case of the following well-known proposition, although since we could 
not find a proof in the literature we include
one here. 

Let $X$ be an affine variety with a rational action of a reductive group $G$. Fix a character 
$\chi :G \longrightarrow \C^\times$  and define the \emph{Zariski open
set of semistable points} to be,
cf. e.g. \cite{Mu}, ch. 6-7, 
$$X^{ss} = \{ x \in X: \text{there exists $m>0$ and $f\in \C[X]^{\chi^m}$ such that $f(x)\neq 0$} \},$$ 
and let $j: X^{ss} \hookrightarrow X$ be the inclusion. By definition the G.I.T. quotient 
$X/\!\!/_{\!\chi}G$ is $\Proj B$ where $B$ is the graded algebra $B = \bigoplus_{m\geq 0} \C[X]^{\chi^m}$. 
Write $\Coh^{G}(Y)$ for the abelian category of $G$-equivariant coherent
sheaves on a $G$-variety $Y$. 

\begin{proposition} \label{restriction-prop}   
{\rm (1)}
For any $\SC\in \Coh^G(X)$,  there exists an integer
$m(\scr S)$ such that  the restriction map induces an isomorphism
\beq{ZZ}
j^*:\
\Gamma(X,\SC)^{\chi^{m}}
\iso\Gamma(X^{ss},\ j^*\SC)^{\chi^{m}}\qquad
\text{for all } m\geq m(\scr S).
\eeq

{\rm (2)} If the orbit map $p: X^{ss} \longrightarrow X/\!\!/_{\!\chi}G$ is a principal $G$-bundle then there is a 
natural isomorphism 
\beq{Z}
j^*\SC  \ \cong \  p^*\ccirc\qgr\big(\oplus_{m\geq 0}\Gamma(X, \SC)^{\chi^{m}}\big)
\eeq
 of functors from $\Coh^{G}(X)$ to $\Coh^{G}(X^{ss})$.
 
\end{proposition}
\noindent

\proof (1) Let $Z = X\setminus X^{ss}$, a closed subvariety of $X$.
 There is an exact sequence \begin{equation} \label{restricttoopen}
  \Gamma_Z(X, \SC) \longrightarrow \Gamma(X, \SC) \longrightarrow \Gamma(X^{ss}, j^*\SC)
   \longrightarrow H^1_Z(X,\SC).\end{equation}
All modules in this sequence are rational $G$-modules, so taking $\chi^m$ semi-invariants 
produces another exact sequence. We claim that both $\Gamma_Z(X, \SC)^{\chi^m}$ 
and $H^1_Z(X, \SC)^{\chi^{m}}$ are zero for $m\gg 0$. Once we have shown this, then it follows that the map
$\Gamma(X, \SC)^{\chi^{m}} \to \Gamma(X^{ss}, j^*\SC)^{\chi^{m}}$ is an isomorphism 
for large enough $m$, thus confirming (1).

To prove the claim, we define  $S = \C[X\times \C]^G$ where $G$ acts on $X\times \C$ by 
$g\cdot (x,\lambda) = (g\cdot x, \chi^{-1}(g)\lambda)$. It follows that $S$ is a finitely generated and so we 
can find a finite set of homogeneous elements $f_1, \ldots , f_n$ in $S_+$ that generate the algebra $S$ 
over $S_0 = \C[X]^G$. Let $I$ be the ideal of $\C[X]$ generated by $f_1, \ldots , f_n$ and observe that $S$ 
has been constructed so that 
 $$Z = \{ x\in X: f(x) = 0 \text{ for all } f\in S_+\} = \{ x\in X : f_i(x) = 0 \text{ for }1\leq i \leq n\}.$$ So to 
 calculate the local cohomology groups $H^i_Z$ it is enough to 
calculate the $H^i_I$. 

Let $R = \C[X]$, $B = \bigoplus_{m\geq 0} \C[X]^{\chi^m}$ and $M = \Gamma(\M , \SC)$. We need to 
calculate the homology of the complex 
$$\begin{array}{rl}
0\longrightarrow M \longrightarrow ({\bigoplus}_{i} R_{{{f}_i}} )\otimes_{R}M \longrightarrow 
 ({\bigoplus}_{i<j} R_{{{f}_i}{{f}_j}} )\otimes_{R}M  & \longrightarrow \cdots \\ \noalign{\vskip 5pt}
  \cdots \longrightarrow  &
 R_{{{f}_1} \cdots {{f}_r}}\otimes_R M\longrightarrow 0.
 \end{array}$$ 
As each $f_i\in B$, we can replace $R$ by $B$ in the above sequence. Since we are going to take 
$\chi$ semi-invariants it is enough for us to study the complex with $M$ replaced by 
$N = \bigoplus_{m\in \Z} M^{\chi^m}$. Set $N^+ = \bigoplus_{m\geq 0} M^{\chi^m}$. Since 
$N/N^+$ is $f_{i_1}\ldots f_{i_r}$-torsion for $r\geq 1$, tensoring the short exact sequence 
$0\longrightarrow N_+ \longrightarrow N \longrightarrow N/N_+ \longrightarrow 0$ by 
$S_{f_{i_1}\ldots f_{i_r}}$ shows that 
$S_{f_{i_1}\ldots f_{i_r}}\otimes_S N_+ \cong  S_{f_{i_1}\ldots f_{i_r}}\otimes_S N$. 
Thus for our calculation we can even replace $N$ by $N_+$ and we need only calculate the 
local cohomology groups $H^i_{S_+}(N_+)^{\chi^m}$.  Now $N_+$ is a finitely generated 
graded $S$-module so thanks to \cite[15.1.5]{BS}, these groups vanish for $m\gg 0$. 
 
 (2) This follows from (1) and the projection formula. In more detail, by part (1)  we have 
$$\qgr\big(\bigoplus_{m\geq 0}\Gamma(X, \SC)^{\chi^{m}}\big) =
 \qgr\big(\bigoplus_{m\geq 0}\Gamma(X^{ss}, j^*\SC)^{\chi^{m}}\big).$$ 
 Since $p$ defines  a principal $G$-bundle we deduce that $p^*[(p_* j^*\SC)^G] \cong j^*\SC$
 (see  \cite[Chapter~VII, Section~1]{SGA1}, again). 
 Thus it is sufficient to check that 
 $$ (p_* \mathcal{F})^G \cong \qgr\big(\bigoplus_{m\geq 0}\Gamma(X^{ss}, \mathcal{F})^{\chi^{m}}\big)
\qquad \text{ for any } \mathcal{F}\in \Coh^G (X^{ss}).$$ But by definition
  $(p_*\mathcal{F})^G \cong \qgr (\bigoplus_{m\geq 0} \Gamma (X/\!\!/_{\!\chi}G , (p_*\mathcal{F})^G(m)).$ 
  Now   \begin{eqnarray*}
   (p_*\mathcal{F})^G(m) \  \cong  \ (p_*\mathcal{F}(m))^G \ 
  & \cong & p_*(\mathcal{F} \otimes_{\mathcal{O}_{X^{ss}}} p^*\mathcal{O}_{X/\!\!/_{\!\chi}G}(m))^G   
   \\   \noalign{\vskip 7pt}
    &&   \cong  \   p_* (\mathcal{F} \otimes \C_{\chi^{-m}})^G \  \cong \ (p_*\mathcal{F})^{\chi^{m}}.
    \end{eqnarray*} 
    The proof is now completed by the following isomorphisms
     \begin{align*}
     \Gamma(X/\!\!/_{\!\chi}G, (p_*\mathcal{F})^G(m)) \cong
       \Gamma(X/\!\!/_{\!\chi}G, (p_*\mathcal{F})^{\chi^{m}}) 
     \cong  \Gamma(X/\!\!/_{\!\chi}G, p_*\mathcal{F})^{\chi^{m}} \cong   
\Gamma(X^{ss}, \mathcal{F})^{\chi^{m}}.\tag*{$\Box$}
     \end{align*}

\subsection{Characteristic cycles}  We assume that $c+m$ is good for all $m\in \mathbb{Z}_{\geq 0}$.
Given a finitely generated left $\Uc$-module $M$, there are two ways to associate to $M$
an algebraic cycle in $\Hilb$.
\medskip

{\bf GS construction} (see \cite[Section~2.7]{GS2} for more details).
Take a   good  filtration  on $M$ and  let $\Lambda$ be the induced 
  tensor product filtration on each 
 $\FM_m =  (\PPP{c+m}{c})\otimes_{\Uc} M$ and hence on 
 $\FM = \bigoplus_{m\geq 0} \FM_m$.   
  By Theorem~\ref{main-thm}, $\gr(\PPP{c+\ell+m}{c+m}) = A^{\ell}$ for all $\ell,m\geq 0$ 
  and it follows easily that $\gr_\Lambda\FM=\bigoplus \gr_\Lambda  \FM_m $ 
  is a graded module  over $A=\bigoplus A^\ell$.   Let $\qgr(\gr\FM)$ denote the 
corresponding coherent sheaf on  $\Proj A=\Hilb$.  We define
$\ch^{GS}(M)$ to be the characteristic cycle of $\qgr(\gr\FM)$. In other words   
$\ch^{GS}(M)$  is the characteristic variety of $\qgr(\gr\FM)$, counting multiplicities, see \cite[(2.7.1)]{GS2}.
 
\medskip
{\bf GG construction} (see \cite[Section 7.5]{GG} for the details).
Recall the definition of  ${}^\top\widetilde{\BH}_c$ from \eqref{commute1}
and  consider the  left 
$\D(\V)$-module 
$${}^\top\widetilde{\BH}_c(M) =( \D_{-c}(\X)/\D_{-c}(\X)\tau(I_{-c}))\otimes_{U_{-c-1}} 
\Omega_c(M).$$ Let $CC({}^\top\widetilde{\BH}_c(M))$ $\sset T^*\X$ be the characteristic cycle of
that $\D_{-c}(\X)$-module,  a closed $PGL(V)$-invariant algebraic cycle
set-theoretically contained in $\mu_{\X}^{-1}(0)$. Following
\cite[Section~7.5]{GG}\footnote{Actually, we have twisted by $\Omega_c$} and in the notation of 
\eqref{hilb-sect}, define 
 $\ch^{GG}(M)$  to be the  unique algebraic  cycle in $\Hilb$ such that 
one has $j_{\X}^* CC({}^\top\widetilde{\BH}_c(M))=p_{\X}^*\ch^{GG}(M)$ inside  $\M\cyc$.
 
\subsection{}  We are now ready to prove the following slight strengthening of
Theorem~\ref{ccthm} from the introduction. 

\begin{theorem} \label{ccthm1} 
Assume that  $c\in \C\smallsetminus  \mathbb{Z}_{<0}$ is chosen such
 that  $c+m$ is good for all integers $m \geq 0$. If $n=2$ 
assume that $nc\not\in \mathbb{Z}$. Then, for any finitely generated  $\Uc$-module $M$, one has an
equality
of algebraic cycles $\ \ch^{GS}(M)=\ch^{GG}(M)$.
\end{theorem}

\proof  Our hypotheses ensure that $nc\not\in \mathbb{Z}_{<0}$, so Corollary~\ref{commute1} 
is available to us.

Put  $F= \DD_{-c}\o_{{\sf U}_{-c-1}}\Omega_{c}(M)$ and note that $F$ is 
$GL(V)$-equivariant since the same is true of  $\DD_{-c}$. For any integer $m\geq 0$ the
 $\D(\V)$-module structure on $F$ induces a natural $\D_{-c-m}(\X)$-module structure
on $F^{(-m)}$.  We compute
\begin{eqnarray*} 
F^{(-m)}&=&[\DD_{-c}\o_{{\sf U}_{-c-1}}\Omega_c(M)]^{(-m)}  \ = \
\DDD{-c-m}{-c}\o_{{\sf U}_{-c-1}}\Omega_c(M) \\  \noalign{\vskip 6pt} &=&
 \DDD{-c-m}{-c}\o_{\D_{-c}(\X)}\left( \frac{\D_{-c}(\X)}{\D_{-c}(\X)\tau(I_{-c})}\right)\o_{{\sf U}
_{-c-1}}\Omega_c(M) \ = \ \BS^m\ccirc\!{}^\top\widetilde{\BH}_c(M),
\end{eqnarray*}
where in the last equality we have used the  analogue of Lemma~\ref{D-routine}(1)  
for   the $\D_{-c}(\X)$-module ${}^\top\widetilde{\BH}_c(M)$.
Hence, taking $SL$-invariants and using Corollary~\ref{commute1}, we deduce that
\begin{align}\label{hdh}
\Omega_{-c-m-1}([F^{(-m)}]^{SL})
&=  \Omega_{-c-m-1}([\BS^m\ccirc{}^\top\widetilde{\BH}_c(M)]^{SL})\\ 
&=  \widetilde{\BH}_{c-m}\ccirc\BS^m\ccirc{}^\top\widetilde{\BH}_c(M) \ = \ 
\PPP{c+m}{c}\o_{\Uc}M.\nonumber
\end{align}

Pick a good filtration  on $M$ and let $\Lambda$ denote the induced tensor product filtration 
on $\PPP{c+m}{c}\o_{\Uc}M$. Writing `$\supp(-)$' for the 
support-cycle, we have
\beq{GS}
\ch^{GS}(M)=\supp\qgr(\gr_\Lambda\FM)\quad\text{where}\quad
\FM =  \oplus_{m\geq
0}\ (\PPP{c+m}{c}\o_{\Uc}M).
\eeq

Since $\phi_c$ preserves the differential operator filtration, the   filtration on $M$ also gives a filtration on 
$\Omega_c(M)$, and hence a  tensor product filtration $\nu$ on $F$.  Observe that, viewed as a $\D(\V)
$-module, this is a good filtration on $F$ since, by \cite[Lemma~6.7(2)]{GS}, $\gr_{\nu}F$ is a 
homomorphic image of $\gr_\Gamma\DD_{-c} \otimes_{\gr U_{c-1}} \gr(\Omega_c(M))$. The filtration on 
$F$ is $GL(V)$-stable and, for any $m\geq 0,$ 
it restricts to a filtration on each of the subspaces
 $[F^{(-m)}]^{SL}\sset F^{(-m)}\sset F$. 
Now, after applying $\Omega_{-c-m-1}$, the composite isomorphism
in \eqref{hdh} transports this filtration on
$\Omega_{-c-m-1}([F^{(-m)}]^{SL})$ to a certain $\nu$-filtration on
$\PPP{c+m}{c}\o_{\Uc}M$. However, this  need not  equal to the 
$\Lambda$-filtration introduced earlier.

Since the action of $GL(V)$ on $F$ is locally
 finite, taking the associated graded $\gr_\nu(-)$  commutes with
taking $GL(V)$-semi-invariants. 
Hence we have
\begin{eqnarray} \label{gradeds} \gr_\nu\FM&=&\bigoplus_{m\geq
0}\ \gr_\nu\left( \PPP{c+m}{c}\o_{\Uc}M\right)=\bigoplus_{m\geq
0}\gr_\nu \left(\Omega_{-c-m-1}( [F^{(-m)}]^{SL})\right) \nonumber \\
&=& \bigoplus_{m\geq
0}\gr \Omega_{-c-m-1} \left(\gr_\nu ( [F^{(-m)}]^{SL})\right) \nonumber 
=\bigoplus_{m\geq
0}[(\,\gr_\nu F)^{(-m)}]^{SL},
\end{eqnarray}
where for the last equality we used that $\gr \Omega_{-c-m-1} $ is the identity map since the associated graded of 
the mapping $z\mapsto \phi_{-c-m-1}(\delta^{-1} z \delta)$ is the identity on 
$\C[\h\times  \h^*]^We$. Since $\gr_\nu F$ is a finitely
generated $\gr_\Gamma \D(\V) = \C[T^*\V]$-module the final object of \eqref{gradeds}
is a finitely generated graded $A$-module. It follows
that the $\nu$-filtration on $\FM$ is good and so
the associated graded modules $\gr_\Lambda\FM$ and $\gr_\nu \FM$ give rise to
the same class in the Grothendieck {\em semi}group of the
category $\Lmof{A}$. We conclude
that
\beq{GS2}
\supp\qgr(\gr_\Lambda\FM)=\supp\qgr(\gr_\nu\FM)=
\supp\qgr\Big(\bigoplus_{m\geq
0}\left[\gr_\nu F^{(-m)}\right]^{SL}\Big).
\eeq

Since $T^*\V$ is affine,
we can write $\gr_\nu F=\G(T^*\V, \CF)$ for a unique
 coherent sheaf $\CF$ supported on the subvariety
 $\M\sset T^*\V$. The sheaf $\CF$ is automatically $GL(V)$-equivariant since $F$ is. 
Taking $SL$-invariants and
applying \eqref{ZZ},  we deduce that,  for large enough $m$, 
$$
\left[\gr_\nu F^{(-m)}\right]^{SL}=
\left[\G(T^*\V, \CF)^{(-m)}\right]^{SL}=
\G(\M,\CF)^{\det^{-m}}=
\G(\M\cyc,\CF)^{\det^{-m}}.
$$
Thus, in the category $\Coh(\Hilb)$, we have the  isomorphisms
$$
\qgr(\gr_\nu\FM)\ \cong\ \qgr\Big({\bigoplus_{m\geq 0} }
\left[\gr_\nu F^{(-m)}\right]^{SL}\Big)\ \cong\
\qgr\Big(\bigoplus_{m\geq
0}\G(\M\cyc,\CF)^{\det^{-m}}\Big).
$$
Applying the isomorphism of functors in \eqref{Z}, we obtain
$$p^*\qgr(\gr_\nu\FM)\ \cong\  p^*\qgr\Big(\bigoplus_{m\geq
0}\G(\M\cyc,\CF)^{\det^{-m}}\Big) \ = \ j^*\CF.
$$

Using this  and \eqref{GS} and \eqref{GS2}, we finally obtain a chain of equalities of
 algebraic cycles
$$
p^*(\ch^{GS}(M)) \ = \ \supp(p^*\qgr(\gr_\nu\FM))\ =\ \supp(j^*\CF).
$$

It remains to prove that $\supp(j^*\CF) \ = \ p^*(\ch^{GG}(M)).$
This is now completely formal. By the commutative diagram \eqref{hilbcomm} and the definition of $\ch^{GG}(M)$
 we have $$p^*(\ch^{GG}(M)) = (\iota\cyc)^* p_{\X}^* (\ch^{GG}(M)) =
  (\iota\cyc)^*j_{\X}^* CC({}^\top\widetilde{\BH}_c(M)) = j_0^*\iota^*  CC({}^\top\widetilde{\BH}_c(M)).$$
Since $\supp(j^*\CF) = j^* \supp(\CF) = j_0^* \supp (j_1^*\CF)$ it is  sufficient to show that 
$\supp (j_1^*\CF) \ = \ \iota^*  CC({}^\top\widetilde{\BH}_c(M)).$ We compute
\begin{align*}
\iota^*  CC({}^\top\widetilde{\BH}_c(M)) =\iota^* \supp ( (\gr_{\nu} F)^{\C^\times}) &=
 \supp ( \iota^*  (\gr_{\nu} F)^{\C^\times} )\\
&=
  \supp ( j_1^* \gr_{\nu} F ) =  \supp ( j_1^* \CF).
\tag*{$\Box$}
\end{align*}

  \section{Appendix: The radial parts map}\label{appendix-section}
  
  \subsection{}  As  was mentioned earlier, in this appendix we will give the details behind 
  Theorem~\ref{first-identity}, since it does not exactly follow from the results in \cite{GG}.
  Thus, in the notation of Section~\ref{sub} our aim is to prove:
  
  \begin{theorem}\label{first-identity1}     The radial components map
  $\Rad = \Rad_c :  \D(\GG\cyc)^G\to\D(\h/W)$ 
  given by 
$$   \Rad(D)=\svol^{c}\ccirc  \bigl(D|_{\mathcal{O}(\GG\cyc,c)}\bigr)\ccirc \svol^{-c}  
\qquad \text{for }  D\in \D(\GG)^G$$
 induces a filtered isomorphism 
   $\Rd : \DD_c^G = \left(\D(\GG)/\D(\GG)\cdot I_c\right)^G \iso \UU_{c-1}.$ 
   
 Consequently, 
   the associated graded map $\gr \Rd: \gr\DD_c^G\to \gr \UU_{c-1}=\C[\h\times\h^*]^W$ is also an isomorphism.
   \end{theorem}

   The proof of this result   closely follows the proof of the analogous results in 
   \cite{EG} and \cite{BFG} (see, in particular \cite[Proposition~5.4.1]{BFG}),
    so the important point here is to determine \emph{which} 
   spherical algebra $\UU_{d}$ contains $\mathrm{Im}(\Rad )$.
   
   To begin, let $\Delta_\g$ be the second order Laplacian on $\g$ associated to  
   non-degenerate invariant bilinear form $(-,-)$ and identify $\Delta_\g$ with $\Delta_\g \otimes 1\in \D
(\GG)$
   acting trivially in the $V$-direction.
   Write  $\Delta_\h$ for the analogous Laplacian on $\h$.
   As usual, we let $\{e_\alpha : \alpha\in R\}$ denote the root vectors for $\g$, normalised so that 
   $(e_\alpha, \,e_{-\alpha})=1$ and set $h_\alpha = [e_\alpha,\, e_{-\alpha}]$ for such $\alpha$.
     Now $\Delta_\g\in \D(\g)^G$ and, as in the analogous computations
   in \cite[Proposition~6.2]{EG},  the 
   key to proving the theorem is to compute  $\Rad ( \Delta_\g)$.

   To do this we slightly change our perspective on $\Rad$.
   Fix a scalar $\cc\in \C$ (which will eventually become $\cc=-c$)
   and let $\GGr = \{ (X,v)\in  \GG^{\cyc} : X\in \g^{\mathrm{rs}}\}$, where $ \g^{\mathrm{rs}}\supset \hreg$ 
  denotes the  regular  semisimple elements.   As in Section~\ref{sub} or 
  \cite[Section~5.4]{BFG} the  projection $\GG\to \g$
   induces  an isomorphism   $\rho^* : \C[\hr/W] \to  \C[\GGr]^G$.    
 For $D\in \D(\GG)^G$ we  define $\Rdd(D)\in \D(\hreg)^W=\D(\hreg/W)$
   to be $\Rad_{-w}(D)$; equivalently, 
  $$\Rdd(D)(f) = \svol^{-\cc}(D(\svol^\cc(\rho^* f)))\big|_{ \hr}
  \qquad \text{for} \ f\in \C[\hreg/W].$$

 \subsection{The radial part of the Laplacian.} 
   We wish to compute $\Rdd(\Delta_\g)$ and we begin by following  the proof of  
   \cite[Proposition~6.2]{EG}.
As observed there, we have the expansion
    \begin{equation}\label{sc0}
     \Delta_{\g} =   \Delta_{\h} \ + \ \sum_{\alpha\in
R}\frac{\partial^2}{\partial e_{\alpha} \partial e_{-\alpha}}
\end{equation}  and  we begin by understanding the final term of this equation.

The differential of the $G$-action on $\GGr$ gives an action of $\g$ and $U(\g)$ on 
$\C[\GGr]$ which we will write as $x\cdot f$ for 
  $x\in U(\g)$ and $f\in \C[\GGr]$. By definition,  $\rho^*f\in \C[\GGr]^G$ for $f\in \C[\hr/W]$ and so   
  $g \defn \svol^\cc \rho^*(f)$ is a $(\det^{-\cc})$-semi-invariant function, since the same is true of 
 $\svol^\cc$. Therefore,  for $X\in \hreg$,  $ue_\alpha \in \C e_\alpha$ and $v\in V$ we have
\begin{eqnarray*} 
(e^{te_{\alpha}}\cdot g)(X+ ue_{-\alpha}, v) & = &
 \bigl(   \det(e^{te_{\alpha}})\bigr)^{-\cc}g(X+ue_{-\alpha},v) \hfill \\
  \noalign{\vskip 5pt}
   & & =  \bigl(1 -  \cc t\Tr(e_\alpha) + O(t^2)\bigr) g(X+ue_{-\alpha},v) ..
\end{eqnarray*}
On the other hand,  
\begin{eqnarray*}  
 &(e^{te_{\alpha}}\cdot g)(X+ ue_{-\alpha},\,  v) \ = \
 g\bigl( (\operatorname{Ad} e^{-te_\alpha})(X+ue_{-\alpha}\bigr),\,  e^{-te_\alpha}\cdot v\bigr)\qquad\qquad \\
   \noalign{\vskip 5pt}
& \qquad =  g\bigl( X + ue_{-\alpha} +t[e_{\alpha}, X]   - tu[e_\alpha,
e_{-\alpha}]  + O(t^2), \, \,  v - te_\alpha\cdot v + O(t^2)\bigr) \\
  \noalign{\vskip 5pt}
& \qquad =  g\bigl( X + ue_{-\alpha} -  t\alpha(X)e_{\alpha}  - tuh_\alpha  + O(t^2), \,\,   v - te_\alpha\cdot v + O
(t^2)\bigr).
\end{eqnarray*}

Now clearly $0=\Tr(e_\alpha) $ and so   after equating the last two equations, applying  $d/dt$,   
and then setting   $t=0$,  we   obtain
\begin{eqnarray*}   
 0  & = &
\alpha(X)\frac{\partial}{\partial e_{\alpha}} g(X+ue_{-\alpha}, v)
- u \frac{\partial}{\partial h_\alpha} g(X+ ue_{-\alpha}, v) -
e_{\alpha}\cdot g(X+ue_{-\alpha}, v) .
\end{eqnarray*}
Rewriting this gives
\begin{equation}\label{sc1}
 \frac{\partial}{\partial  e_\alpha}g(X+ue_{-\alpha}, v)  \ =\   \frac{1}{\alpha(X)}\left(
u\frac{\partial}{\partial h_{\alpha}} g(X+ue_{-\alpha},v) +
e_\alpha \cdot g(X+ue_{-\alpha},v)\right).
\end{equation} 
Now   apply  $d/du$  to this equation and then set $u=0$ to  give 
\begin{eqnarray}\label{sc2}
\frac{\partial^2}{\partial e_{-\alpha}\partial e_{\alpha}} g(X,v)
\ = \  \frac{1}{\alpha(X)} \left[  \frac{\partial }{\partial h_\alpha}
g(X,v) + \frac{\partial}{\partial e_{-\alpha}} (e_\alpha\cdot
g(X,v))\right]. \end{eqnarray}
Applying  \eqref{sc1}  with  $e_\alpha$  replaced  by   $e_{-\alpha}$
gives 
\begin{eqnarray*}
 \frac{\partial }{\partial e_{-\alpha}} g(X, v+
te_\alpha \cdot v) & = & -\frac{1}{\alpha(X)} e_{-\alpha}\cdot g(X,
v+ te_\alpha \cdot v).
 \end{eqnarray*}
 Therefore \ 
$\displaystyle \frac{\partial}{\partial e_{-\alpha}} (e_\alpha \cdot g(X,v) ) \ = \ 
-\frac{1}{\alpha(X)} e_{-\alpha}e_\alpha \cdot g(X,v)$ and so \eqref{sc2} becomes
 \begin{eqnarray}   \label{sc3}
  \frac{\partial^2}{\partial  e_{-\alpha}\partial e_{\alpha}} g(X,v)  \ = \  \frac{1}{\alpha(X)}
\frac{\partial }{\partial h_\alpha} g(X,v) - \frac{1}{\alpha(X)^2}
e_{-\alpha}e_\alpha\cdot g(X,v)\end{eqnarray}
 for $X\in \hr$ and $v\in V$.

 Next we calculate $e_{-\alpha}e_\alpha \cdot g(X,v)$ for 
$X\in \hr$. Recall  that $g = \svol^\cc \rho^*(f)$ with 
$\rho^*f\in \C[\GGr]^G$. As noted in \cite[Section~5.4]{BFG},
projection onto the first term gives an isomorphism $\C [\GG]^G\ \buildrel{\sim}\over{\too}\ \C[\g]^G$
and so   $e_{-\alpha}\cdot \rho^*(f)\big|_{\hreg} = e_\alpha \cdot
\rho^*(f)\big|_{ \hr} = 0$.  Thus 
$$e_{-\alpha}e_{\alpha}\cdot
\svol^\cc \rho^*(f)\big|_{ \hr}  \ = \ (e_{-\alpha}e_{\alpha}\cdot \svol^\cc \big|_{
\hr}) f.$$ 
Now for $X=(x_1,\dots, x_n)  \in \hr$  and $v=(v_1,\dots,v_n)\in V$
we have \begin{equation} \label{sc}
\svol^\cc  (X,v) \ =\  \prod_{i < j} (x_i-x_j)^\cc  (v_1\cdots
v_n)^\cc .\end{equation} 
If   $\alpha $ is the elementary matrix $\alpha= E_{ij}$ then 
\begin{eqnarray*}e_\alpha \cdot (v_1^{a_1} \cdots v_n^{a_n}) &=&
   a_i v_1^{a_1} \cdots  v_i^{a_i-1}  \cdots
v_j^{a_j+1} \cdots v_n^{a_n} \qquad\text{for any } a_\ell \in \mathbb{C}.
\end{eqnarray*}
Thus  $e_{-\alpha}e_{\alpha}\cdot \svol^w\big|_{ \hr} = w(w+1)\svol^w\big|_{
\hr} $   and so
\begin{equation}\label{sc4}
e_{-\alpha}e_{\alpha}\cdot
g(X,v)\big|_{ \hr} \ = \  \cc (\cc +1)g(X,v)\big|_{ \hr}.
\end{equation}
Finally, observe that   $\svol^\cc \big|_{ \hr} = \delta^\cc $. Thus combining \eqref{sc0},
\eqref{sc3} and \eqref{sc4}  produces the desired equation
\begin{equation}\label{sc5} 
\Rdd(\Delta_\g) \ =\  \delta^{-\cc } \Delta_{\h} \delta^\cc  \  +\ \sum_{\alpha\in R}
\frac{\delta^{-\cc }}{\alpha}\frac{\partial}{\partial h_\alpha}
\delta^\cc  \ - \  \sum_{\alpha \in R} \frac{\cc (\cc +1)}{\alpha^2}. 
\end{equation}

\subsection{}\label{onto-subsect}  Recall from \eqref{cherednik-defn} that we have identified $\HH_{\cc}$ 
with its image in $ \D(\hreg)\ast W$ under  the Dunkl  embedding.
 As in  \cite[p.~281]{EG}, we define a  \emph{twisted Dunkl homomorphism}  
$\Theta_\cc ^{\spher} : \HH_{\cc}\to  \D(\hreg)\ast W$ by 
$\Theta_\cc ^{\spher}(h) = \delta^{-\cc } h  \delta^{\cc }$ for $h\in \HH_\cc$.

Recall further the  \emph{Calogero-Moser
operator} defined by the formula
$$L_\cc  \  \defn\  \Delta_\h - \frac{1}{2}\sum_{\alpha\in R} \cc (\cc +1) \frac{(\alpha ,
\alpha)}{\alpha^2}.$$ 

It was mentioned in \cite[pp.~281-2]{EG}  that  the  Calogero-Moser
operator lies in  $\Theta^{\spher}_\cc  (\HH_\cc ).$
 We take this opportunity to
give a complete proof of this assertion. In this computation  we write  $\partial_t = \partial/\partial{x_t}$
for $t\geq 1$ and $x_{ij}=x_i-x_j$ for $i\not=j$. Then
  \eqref{dunkdef} gives 
 \begin{eqnarray*}
 \sum_{t=1}^n D_w(y_t)^2 e &=&
  \sum_{t=1}^n \bigg( (\partial_t)^2 - \sum_{j<i} w\frac{\langle y_t, \, x_{ij} \rangle}{x_{ij}} (\partial_t - s_{ij}(\partial_t))\bigg) \\ 
  &=& \sum_{t=1}^n \bigg( (\partial_t)^2 - 2 \sum_{j<i} w\frac{\langle y_t, \, x_{ij} \rangle}{x_{ij}} \partial_t\bigg).
  \end{eqnarray*} 
  \medskip 
  
  Now 
 $\delta^{-w}\sum_{t=1}^n \partial_t = \sum_{t=1}^n \bigg( \partial_t \delta^{-w} +
  \sum_{j<i} w \frac{\displaystyle \langle y_t, \, x_{ij} \rangle}{
 \displaystyle x_{ij}} \delta^{-w}\bigg).$ 
Thus, we compute
  \allowdisplaybreaks\begin{align*}
\delta^{-w} &\sum_{t=1}^n \bigg( (\partial_t)^2 - 2 \sum_{j<i} 
w\frac{\langle y_t, \, x_{ij} \rangle}{x_{ij}} \partial_t\bigg)\\
&=    \sum_{t=1}^n \bigg( \partial_t\delta^{-w}\partial_t  + \sum_{j<i} w
 \frac{\langle y_t, \, x_{ij}  \rangle}{x_{ij}}
     \delta^{-w}\partial_t - \delta^{-w}2 
\sum_{j<i} w\frac{\langle y_t, \, x_{ij} \rangle}{x_{ij}} \partial_t\bigg) \\
 \noalign{\vskip7pt} &=
   \sum_{t=1}^n \bigg( \partial_t\delta^{-w}\partial_t  - \sum_{j<i} w \frac{\langle y_t, \, x_{ij} \rangle}{x_{ij}}
    \delta^{-w}\partial_t \bigg) \\
 &=   \sum_{t=1}^n \bigg( \partial_t^2\delta^{-w} + \sum_{j<i} w \partial_t \frac{\langle y_t, \, x_{ij} \rangle}{x_{ij}}
     \delta^{-w} -  \sum_{j<i} w \frac{\langle y_t, \, x_{ij} \rangle}{x_{ij}} \delta^{-w}\partial_t \bigg) \\
&=  \sum_{t=1}^n \bigg( \partial_t^2\delta^{-w} + \sum_{j<i} w \frac{\langle y_t, \, x_{ij} \rangle}{x_{ij}}\partial_t \delta^{-w} 
   - \sum_{j<i} w \frac{\langle y_t, \, x_{ij}\rangle^2}{x_{ij}^2}   \delta^{-w} 
          -  \sum_{j<i} w \frac{\langle y_t, \, x_{ij} \rangle}{x_{ij}} \delta^{-w}\partial_t \bigg)   \\
&=       \sum_{t=1}^n \bigg( \partial_t^2\delta^{-w} + \sum_{j<i} w \frac{\langle y_t, \, x_{ij} \rangle}{x_{ij}}[\partial_t, \delta^{-w}]
        - \sum_{j<i} w \frac{\langle y_t, \, x_{ij}\rangle^2}{x_{ij}^2}  \delta^{-w} \bigg) \\
&=         \sum_{t=1}^n \bigg( \partial_t^2\delta^{-w} - \bigg(\sum_{j<i} w \frac{\langle y_t, \, x_{ij} \rangle}{x_{ij}}\bigg)^2\delta^{-w}
         - \sum_{j<i} w \frac{\langle y_t, \, x_{ij}\rangle^2}{x_{ij}^2}   \delta^{-w} \bigg) \\ 
 \noalign{\vskip7pt} &=
         \sum_{t=1}^n \bigg( \partial_t^2\delta^{-w} - \!\!\!\!\!\!\!\sum_{\substack{  i< j, k< \ell \\ (i,j)\neq (k,\ell)}}
          \!\!\!\!\!\!\!w^2\frac{\langle y_t, \, x_{ij}\rangle
          \langle y_t, x_{k\ell}\rangle}{x_{ij}x_{k\ell}}\delta^{-w} 
                 -   \sum_{j<i} (w^2+w) \frac{\langle y_t, 
\, x_{ij}\rangle^2}{x_{ij}^2}   \delta^{-w} \bigg) \\ 
 &=      \bigg( \sum_{t=1}^n \delta_t^2 -\!\!\!\!\!\!\!\sum_{\substack{i< j, k< \ell \\ (i,j)\neq (k,\ell)}}
      \!\!\!\!\!\!\!w^2\frac{(  x_{ij},\,x_{k\ell})}{x_{ij}x_{k\ell}} 
      - \sum_{j<i} (w^2+w) \frac{( x_{ij}, \, x_{ij})}{x_{ij}^2} \bigg)\delta^{-w} \\ &=
       \bigg(\Delta_{\h} - \frac{1}{2}(w^2+w) \sum_{\alpha\in R} \frac{( \alpha , \alpha ) }{\alpha^2}\bigg)\delta^{-w}.
\end{align*}
Here the middle term of the second-to-last line disappears since if we write  the sum over the common denominator 
$\delta$ then the numerator becomes a sign semi-invariant element in $\C[\h]$ of degree $\deg(\delta) - 2$, and so is
zero. We deduce that
$ L_w = \Theta^{\spher}_\cc( \sum_{t=1}^n D_w(y_t)^2e),$
as claimed. 

For Cherednik algebras of  type $A$ one has   $(\alpha, \alpha) = 2$ and an easy computation shows 
that 
$\delta \Delta_{\h} \delta^{-1}  = \Delta_\h- \sum_{\alpha\in R}
\frac{\delta}{\alpha}\frac{\partial}{\partial h_\alpha} \delta^{-1} $.  Combined with \eqref{sc5} this shows 
that 
 $$ \delta^{\cc +1} \Rdd(\Delta_\g ) \delta^{-(\cc +1)}\ =\   L_\cc \  =\ L_{-(\cc +1)} \  
    \in \  \Theta_{-(\cc +1)}^{\spher}({\mathsf H}_{-(\cc +1)}).$$
Since $\Rdd(\Delta_g)$ is clearly $W$-invariant  we have therefore proved  that
\begin{equation}\label{8.3.1}
 \Rdd(\Delta_\g)  \in {\mathsf H}_{-(\cc +1)}^W \ = \ \UU_{-(\cc+1)}.
 \end{equation} 

\subsection{} \label{ker-subsect}
To calculate (part of) the kernel of  $\Rdd$  we just
have to see what happens to the element  $\tau(\eu)\in \de$.  As was observed in Section~\ref{hamilton-defn2}, 
$\tau(\eu) = - \sum_{i=1}^n e_i\frac{\partial}{\partial e_i},$ where $\{e_i\} $ is a basis of 
$V^*\subset \C[V]$.  It therefore follows from \eqref{sc} that $\tau(\eu)\cdot\svol^w = -nw \svol^w.$ 
As  before, $\tau(\eu) \rho^*(f)=0 $ for $f\in \C[\hreg/W]$
and so $\tau(\eu)\cdot(\svol^w \rho^*(f)) = -nw(\svol^w\rho^*(f)$. Thus   
\begin{equation}\label{ker-equ}
\mathrm{Ker}(\Rdd)  \ \supseteq \ \bigl(\D(\GG)\cdot\tau(\mathfrak{sl}_n)\bigr)^G 
+ \D(\GG)^G(\tau(\eu )+w\Tr(\eu)).
\end{equation} The right hand side of  \eqref{ker-equ}  is 
   the ideal  $\bigl(\D(\GG)I_{-w}\bigr)^G$ in the notation of  Section~\ref{2.2} and so 
   $\Rad$ does at least induce a homomorphism $\Rd: \DD^G_{-w}\to \ureg$. 
 
\subsection{Proof of Theorem~\ref{first-identity1}}
The proof now follows  a well-worn path similar to that  of   \cite[Theorem~4.8 and Proposition~4.9]{EG} 
or   \cite[Proposition~5.4.4]{BFG}.  Set $c=-\cc$ and $\Rad=\Rad_c$.   Then   \eqref{8.3.1} shows
 that $\Rad(\Delta_\g)\in \UU_{c-1}$.  By  the proof of  \cite[Proposition~5.4.4]{BFG} this is enough to 
 ensure that  $\im\Rad =\UU_{c-1}$.  In fact that proof shows rather more: as $\Rad$ is defined by restriction of 
 differential operators 
 it is certainly a filtered morphism and the proof of surjectivity in \cite{BFG} is obtained by proving that 
 the associated graded map $\gr \Rad$ is surjective.  Therefore, $\Rad$ is filtered surjective and hence,
 by \eqref{ker-equ}, so is $\Rd$.
 
 It remains to prove that  $\Rd$ and $\gr \Rd$ are injective.  However, just as in the 
 proof  of Lemma~\ref{sign-reps},
  $\gr \DD_c^G  \cong \bigl(\gr( \D(\GG)/\D(\GG)I_{c})\bigr)^G$ is a homomorphic image of 
 $$\left[ \frac{\displaystyle\text{gr}\D({\GG}) }{\strut\displaystyle\text{gr}\D({\GG})\cdot
\text{gr}\, \tau(I_{c})} \right]^G
\ =\ 
  \left[\frac{\displaystyle \C[T^*{\GG}] }{
\strut \displaystyle\C[T^*{\GG}]\cd \mu_{\GG}^*(\g)} \right]^G 
\ = \   \C[\mathcal{M}]^G  \ = \  \C[\h\times\h^*]^W,$$
 where $\mathcal{M}$ is defined in \eqref{M} and   the final equality is \cite[Lemma~2.11]{GG}.  
 In other words,  $\gr \Rd$ is a surjective morphism from a factor of 
 $\C[\h\times\h^*]^W$ to $\gr \UU_{c-1} = \C[\h\times\h^*]^W$.  Thus  it must be an isomorphism. 
 This also forces $\Rd$ to be an isomorphism.
\qed

  \setcounter{equation}{0}
\footnotesize{
}
\end{document}